\documentclass{article}
\usepackage{amssymb,amsmath}
\usepackage{anysize}
\newtheorem{thm}{Thm.}[section]
\newtheorem{def.}{Def.}[section]

\newtheorem{prop}{Prop.}[section]
\newtheorem{cor}{Cor.}[section]

\input xy.sty
\xyoption{all}
\xyoption{graph}
\xyoption{poly}
\xyoption{knot}
\xyoption{rotate}

\begin{document}

\def\TC{\xygraph{!{0;/r.75pc/:}
!P8"c"{~:{(3.5,0):}~>{}}
!{\ar @{->} "c3";"c6",^r}}}
\def\TD{\xygraph{!{0;/r.75pc/:}
!P4"a"{~:{(1.75,0):}~>{}}!P8"b"{~:{(3.5,0):}~>{}}
!{\hunder~{"a2"}{"a1"}{"a3"}{"a4"}}
!{\ar @{<-} "a2";"b3",^{r'_1}}
!{\ar @{->} "a3";"b6",^{r'_2}}
!{\save0;"a3"-"a4":"a4",
\hcap-@(+.1)\restore}}}
\def\TA1{\xygraph{!{0;/r.75pc/:}
!P8"c"{~:{(3.5,0):}~>{}}
!{\ar @{->} "c3";"c6",^r}}}
\def\TB1{\xygraph{!{0;/r.75pc/:}
!P4"a"{~:{(1.75,0):}~>{}}!P8"b"{~:{(3.5,0):}~>{}}
!{\hover~{"a2"}{"a1"}{"a3"}{"a4"}}
!{\ar @{<-} "a2";"b3",^{r'_1}}
!{\ar @{->} "a3";"b6",^{r'_2}}
!{\save0;"a3"-"a4":"a4",
\hcap-@(+.1)\restore}}}
\def\RA{\xygraph{!{0;/r.6pc/:}
!P4"a"{~:{(-.75,-.75):}~>{}}!P8"b"{~:{(6,0):}~>{}}
!{\POS"b3" \ar @/^3ex/ "b6", _{r}}
!{\POS"b2" \ar @/_3ex/ "b7", ^{s}}}}
\def\RB{\xygraph{!{0;/r.6pc/:}
!P4"a"{~:{(-.5,-.5):}~>{}}!P8"b"{~:{(6,0):}~>{}}
!{\xunderv~{"b2"}{"b3"}{"a2"}{"a4"}@(.65)<>|><{r'_1}|{s'}}
!{\xunderv~{"a4"}{"a2"}{"b6"}{"b7"}@(.35)<{r'_2}>{r'_3}}}}
\def\QA{\xygraph{!{0;/r.6pc/:}
!P4"a"{~:{(-.75,-.75):}~>{}}!P8"b"{~:{(6,0):}~>{}}
!{\POS"b3" \ar @/^3ex/ "b6", _{r}}
!{\POS"b7" \ar @/^3ex/ "b2", _{s}}}}
\def\QB{\xygraph{!{0;/r.6pc/:}
!P4"a"{~:{(-.5,-.5):}~>{}}!P8"b"{~:{(6,0):}~>{}}
!{\xunderv~{"b2"}{"b3"}{"a2"}{"a4"}@(.65)<>|<<{r'_1}|{s'}}
!{\xunderv~{"a4"}{"a2"}{"b6"}{"b7"}@(.35)<{r'_2}>{r'_3}}}}
\def\Reid{
\xy 0;/r1pc/:
\xoverv~{(-1,3)}{(1,3)}{(-1,1)}{(1,1)}|<<{r'}|{s'}>{u'}  \xoverv~{(-4,-2)}{(-4,-.5)}{(0,-.5)}{(-1,1)}<{v'}|{t'}
\xunderv~{(4,-.5)}{(4,-2)}{(1,1)}{(0,-.5)}|<<{w'}
\endxy}

\def\dieR{
\xy 0;/r-1pc/:
\xoverv~{(-1,3)}{(1,3)}{(-1,1)}{(1,1)}<{w}|{v}>{u}  \xoverv~{(-4,-2)}{(-4,-.5)}{(0,-.5)}{(-1,1)}<<|<<{s}|{t}
\xunderv~{(4,-.5)}{(4,-2)}{(1,1)}{(0,-.5)}<{r}
\endxy}

\title{Quandles at Finite Temperatures II}
		\author{F. Miguel Dion\'\i sio\\
		Departamento de Matem\' atica and Centro de L\' ogica e Computa\c c\~ao\\
	        Instituto Superior T\'ecnico\\
	        Av. Rovisco Pais\\
		1049-001 Lisboa\\
	        Portugal\\
	        \texttt{fmd@math.ist.utl.pt}
\and 		Pedro Lopes\\
		Departamento de Matem\'atica and Centro de Matem\'atica Aplicada\\
	        Instituto Superior T\'ecnico\\
	        Av. Rovisco Pais\\
		1049-001 Lisboa\\
	        Portugal\\
		\texttt{pelopes@math.ist.utl.pt}\\
\and 		To Jorge, in memoriam.}

\date{May 6, 2002}
\maketitle

\begin{abstract}
The number of colorings of a knot diagram by a quandle has been shown to be a knot invariant by CJKLS using quandle cohomology methods. In a previous paper by the second named author, the CJKLS invariant was refined and, in particular, it was shown that the number of colorings is an invariant directly without resorting to quandle cohomology. Here we investigate the effectiveness of this invariant on (classical) knots of up to and including ten crossings.
\end{abstract}

\section{Introduction} \label{S:intro}
In order to clarify the proofs of the invariance of the state-sums set forth in Carter et al (CJKLS) \cite{jsCetal} (theorems 4.4 and 5.6) a careful analysis of the notion of coloring of a diagram was developed in \cite{pLopes} and, in particular, an equivalence relation between colored diagrams was introduced. In short, two colorings are said to be equivalent if the underlying diagrams are related by Reidemeister moves, see \cite{lhKauffman}, (resp., Roseman moves, see \cite{dRoseman1}) and at each move the changes to the previous coloring are consistently performed, i.e., respecting the ``coloring condition'', see \cite{jsCetal}. It is then readily seen that given two equivalent diagrams and fixing a finite set of Reidemeister moves (resp., Roseman moves) connecting them, there is a 1-1 correspondence between colorings for the two diagrams. In particular, the number of (different) colorings by a fixed quandle is an invariant for knots (resp., knotted surfaces). 

In this paper we investigate the effectiveness of this invariant on (classical) knots (i.e., embeddings of $S^1$ in $\mathbb{R}^3$) of up to and including ten crossings. The invariant should be regarded, for each knot, as a finite sequence of numbers of colorings, one per quandle of a fixed list. In \cite{FR} it is suggested that one should look for a quandle that colors only trivially one of the knots and has at least one non-trivial coloring for the other knot. As far as we know this project has never been implemented. Also, based on results referred to above, our approach allows us to consider different numbers of colorings as a valid way of telling knots apart. Moreover, we can apply these ideas to knotted surfaces, too (paper in preparation).

We use two different ways of encoding the knots in order to cross-check our results. In one of them, knots are regarded as closures of braids (see \cite{Birman}) and we follow the approach in \cite{jsCetal1}. In this way knots are represented by Burau matrices (see \cite{Birman}) and the solutions to a certain eigenvector equation yield the colorings for each quandle. In the other one, knots are encoded by Gauss codes (see \cite{CetS}). From each Gauss code we set up a system of equations whose solutions yield the colorings. We detail this approach below. Also, we list the Gauss codes used in Appendix B. Finally, we set up a matrix indexed by knots whose entries are the first quandle to have distinguished the corresponding pair of knots (see Appendix A). Of course, we only use the part of the matrix strictly above the main diagonal. The algorithms were implemented with C.

\subsection{Organization and Acknowledgements} \label{S:org}

In Section 2 we recall the definitions of quandle, of coloring of a knot diagram, of relation amongst colorings and some related facts. In Section 3 we describe the methods and algorithms involved in the calculation of the invariant. In Appendix A we present the invariant as entries of a matrix indexed by (classical) knots of up to and including ten crossings. Each such entry is the first quandle yielding different numbers of colorings for the corresponding knots. Due to the size of the matrix, the relevant part of it has been split into submatrices. In Appendix B we list Gauss codes used for the knots as we weren't able to find such a list in the literature.
 
This work was supported by the programme {\em Programa Operacional ``Ci\^{e}ncia, Tecnologia, Inova\c{c}\~{a}o''} (POCTI) of the {\em Funda\c{c}\~{a}o para a Ci\^{e}ncia e a Tecnologia} (FCT), cofinanced by the European Community fund FEDER.

The authors would like to thank Prof. R. F. Picken for his valuable remarks and J. C. T. Bojarczuk for pointing out the benefits of a low level programming language.

\section{Preliminaries} \label{S:pre}
In this section we first recall the definition of a quandle. We then recall the notion of coloring of a diagram and the (equivalence) relation between the colorings of the diagrams of a knot and the fact that the list of equivalence classes of colorings of a diagram is a knot invariant.

\subsection{ The Quandles} \label{S:qdles}

\begin{def.}[Quandle]A quandle is a set $X$ endowed with a binary operation, denoted $\triangleright$, such that:

(i) for any $a \in X$, $a\triangleright a = a$;

(ii) for any $a$ and $b \in X$, there is a unique $x \in X$ such that $a = x\triangleright b$;

(iii) for any $a$, $b$, and $c \in X$, $(a\triangleright b)\triangleright c = (a\triangleright c)\triangleright (b\triangleright c)$ (self-distributivity).
\end{def.}

Quandles have been studied before (see \cite{dJoyce}, \cite{sMatveev}, \cite{lhKauffman}, and \cite{FR}, for instance). The immediate example is any group with group conjugation as the quandle operation. A first class of quandles are those whose underlying set is ${\mathbb Z}_n$ and quandle operation is $a\triangleright b:=2b-a$, (mod $n$); they are denoted $R_n$. Another important class of quandles stems from the fact that the set of Laurent polynomials $\Lambda = {\mathbb Z}[T, T^{-1}]$ endowed with the operation $a\triangleright b = Ta + (1-T)b$ is a quandle. It follows that any module over $\Lambda $ is a quandle. We will be dealing with the so-called finite Alexander quandles which have the form ${\mathbb Z}_n[T, T^{-1}]/(h)$, where $n$ is an integer, and $h$ is a monic polynomial. We remark that the $R_n$ quandles can be regarded as finite Alexander quandles. There is also a related algebraic structure, known as rack, whose binary operation only has to satisfy properties (ii) and (iii) above.

\subsection{The Colorings of Classical Knots} \label{S:cknots}

In this subsection we present a brief description of colorings and the relation between them. For a more complete discussion of this material see \cite{pLopes}.
For the purpose of the discussion in this subsection we will arbitrarily fix a finite quandle $X$, called the ``labelling quandle'', and a knot $K$ (assumed to be oriented); $D$ and $D'$ will stand for diagrams of $K$, identified up to planar isotopy, and $R_D$ and $R_{D'}$ for their sets of arcs.

\begin{def.}[Coloring of a Diagram]A coloring of a diagram, $D$, is a map ${\cal C}:R_D \rightarrow X$ such that, at each crossing:

$$
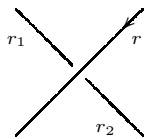
\xy 0;/r1pc/:
,{\xoverv[4.0] |> <{r_1}|{r}>{r_2}}
\endxy$$

if ${\cal C}(r_1)=x$ and ${\cal C}(r)=y$, then ${\cal C}(r_2)=x\triangleright y$. 
%Equivalently, at each crossing:
%
%$$\xy 0;/r1pc/:
%,{\xoverv[4.0] |< <{r_1}|{r}>{r_2}}
%\endxy$$ 
%
%
%if ${\cal C}(r_1)=x$ and ${\cal C}(r)=y$, then ${\cal C}(r_2)=x{\triangleright ^{-1}}y$
%\vskip 5pt

\end{def.}

\vskip 5pt

\begin{def.}[Set of Colorings of a Knot]We will refer to the set of all colorings of all diagrams of a knot simply as the {\bf set of colorings of a knot}.
\end{def.}
\vskip 5pt

Note that, although there may be lots of different ways of coloring two different diagrams, some of these different colorings are related in a natural way:

%RELATION ON THE SET OF COLORINGS OF K

\begin{def.}[Relation $\sim$ on the set of colorings of $K$]Let ${\cal C}:R_D \rightarrow X$ and ${\cal C'}:R_{D'} \rightarrow X$ be colorings. We say ${\cal C}\sim {\cal C'}$ if there is a finite sequence of steps taking ${\cal C}$ to ${\cal C'}$ as described below (using the notation $D$ and $D'$ also to describe each step):
\begin{itemize}

\item If $D'$ is obtained from $D$ by a type I Reidemeister move, then ${\cal C}$ and ${\cal C'}$ are related as follows:
$$\TC\longleftrightarrow\TD$$
with

\[{\cal C}(r)={\cal C'}(r'_1)={\cal C'}(r'_2)
\]
The prescription is analogous for all other cases of type I moves for oriented diagrams. 
%Note that, due to the definition of coloring, the two arcs in the figure on the right have to be assigned the same color, ${\cal C'}(r'_1)={\cal C'}(r'_2)$, for any coloring.

\vskip 50pt

\item If $D'$ is obtained from $D$ by a type II Reidemeister move, then ${\cal C}$ and ${\cal C'}$ are related as follows:
\vskip 10pt
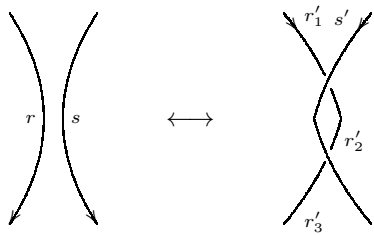
$$\RA\longleftrightarrow\RB$$
\vskip 10pt
\begin{center}
	${\cal C}(r)=x$ \quad and \quad ${\cal C}(s)=y$ \\
	$\sim $\\
	${\cal C'}(r'_1)=x$,\quad ${\cal C'}(s')=y$,\quad ${\cal C'}(r'_2)=x\triangleright y$\quad and\quad ${\cal C'}(r'_3){\triangleright }y=(x\triangleright y)$ 
\end{center}
Analogously for the other cases.
%Another case is:
%\vskip 10pt
%$$\QA\longleftrightarrow\QB$$
%\vskip 10pt
%\begin{center}
%	${\cal C}(r)=x$\quad and\quad ${\cal C}(s)=y$\\
%	$\sim $\\
%	${\cal C'}(r'_1)=x$,\quad ${\cal C'}(s')=y$,\quad ${\cal C'}(r'_2)=x{\triangleright }^{-1}y$\quad and\quad${\cal C'}(r'_3)=(x{\triangleright }^{-1}y)\triangleright y=x$
%\end{center}
%Note that the relation between the colorings of the arcs on the right hand side hold for any coloring ${\cal C'}$, from def. 3.1 and thus the passage from ${\cal C'}$ to ${\cal C}$ is well-defined. The remaining cases of Reidemeister II moves between oriented diagrams are dealt with analogously. 

\vskip 50pt

 \item If $D'$ is obtained from $D$ by a type III Reidemeister move, then ${\cal C}$ and ${\cal C'}$ are related as follows:
%For the third Reidemeister move, let us first remark that each line segment can have two possible orientations, giving a total of $2^3=8$ possibilities for each (unoriented) configuration of the three line segments. Since the colorings do not depend on the orientation of the lower line segment then we just have to take care of $2^2=4$ possibilities. Now for each (unoriented) configuration of the lower line segment there are two possibilities for the remaining arcs (one is the top line segment and the other the middle line segment). Rotation through angles of $\frac{2}{3}\pi$ and $\frac{4}{3}\pi$ gives the remaining (unoriented) configurations; but as we have seen before, this is the same as assigning first colors to other arcs. Thus, a priori, there are eight cases to be specified. Further inspection of these eight cases cuts it down to four. We deal here with one of them remarking that the remaining ones are dealt with analogously.
%\vskip 500pt

$$\dieR
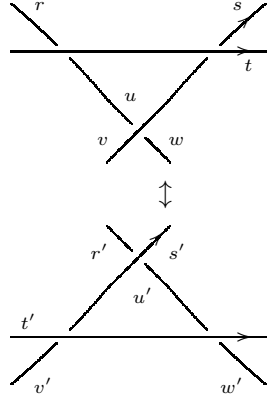
$$
$$\qquad\updownarrow $$
$$\Reid$$

\begin{center}
	${\cal C}(r)=x$,\quad ${\cal C}(s)=y$,\quad and\quad ${\cal C}(t)=z$\quad with\quad ${\cal C}(u){\triangleright }z=x$,\quad ${\cal C}(v){\triangleright }z=y$\quad and\quad ${\cal C}(w){\triangleright }{\cal C}(v)={\cal C}(u)$\\
	$\sim $\\
	${\cal C'}(r')=x$,\quad ${\cal C'}(s')=y$,\quad and \quad${\cal C'}(t')=z$,\quad ${\cal C'}(u'){\triangleright }y=x$,\quad with\quad ${\cal C'}(v'){\triangleright }z=y$,\quad and \quad ${\cal C'}(w') {\triangleright } z={\cal C'}(u')$
\end{center}

Analogously for the other cases.

\end{itemize}

\end{def.}

\vskip 20pt

We remark that:

\begin{prop}$\sim$ is an equivalence relation on the set of colorings of $K$.\end{prop}
%$\hfill \blacksquare$

Suppose $D$ and $D'$ are two (Reidemeister) equivalent diagrams and fix $n$ Reidemeister moves connecting them. Consider the sets $S=\{ {\cal C}: R_D \mapsto X \}$ and ${S'}=\{ {\cal C'}: R_{D'} \mapsto X \}$ of all colorings from the set of arcs of each diagram to a fixed quandle $X$. Consider also the maps:
\[
\begin{array}{clcr}	
	\psi : & S \longrightarrow S'\\
		&{\cal C} \longmapsto {\cal C'}
\end{array}
\quad\txt{ and }\quad
\begin{array}{clcr}
	\psi ': & S' \longrightarrow S\\
		&{\cal C'} \longmapsto {\cal C}
\end{array}
\]
(where, in both cases, ${\cal C} \sim {\cal C'}$ via the above fixed 
Reidemeister moves). It is readily seen that both maps are surjective and 
their composite is the identity map. Hence both of them are bijections 
relating pairs of equivalent colorings and they establish a 1-1 correspondence 
between $S$ and $S'$.

In this way, we have just proved that:

\begin{thm}The set of equivalence classes of the colorings associated to a diagram (along with their multiplicities) is a knot invariant.\end{thm}

\vskip 10pt

\begin{cor}The number of equivalence classes of the colorings associated to a diagram is a knot invariant.\end{cor}

\vskip 10pt

\begin{cor}The set of multiplicities referred to in Thm. 2.1 is a knot invariant.\end{cor}

\vskip 10pt

\begin{cor}The total number of colorings associated to a knot diagram is a knot invariant.\end{cor}

\vskip 10pt

\section{Methods and Algorithms} \label{S:algos}

Based on the results referred to above, we fix a finite sequence of quandles and, for each (unordered) pair of knots we look for the first quandle yielding different numbers of colorings for these knots. The number of this quandle in the sequence then becomes the entry of a matrix indexed by the two knots. We repeat this procedure for all pairs of (unordered) knots of up to and including ten crossings (see \cite{R} and \cite{J}). Note that only the part of the matrix strictly above the main diagonal is needed.

We use two different ways of encoding the knots. In one of them they are regarded as closures of braids and are thus encoded as the Burau matrices representing the corresponding braids (see \cite{Birman}). The Burau matrices represent also the effect of a generic Alexander quandle on the colors introduced on the top of the corresponding braid by multiplying a vector whose coordinates are these colors. These colors stand for a coloring if the Burau matrix leaves them invariant. In this way colorings are given by solutions to an eigenvector equation (and vice versa). The different (Alexander) quandles are obtained by modding out by different monic polynomials and by modding out the polynomial coefficients by different integers. This approach has already been used in \cite{jsCetal1} although the goal there was the computation of the CJKLS invariant.

In the other approach we use Gauss codes to encode the knots. A Gauss code is obtained from a planar diagram of a knot as follows. Fix an arbitrary starting point on an arc of the (oriented) planar diagram. Crossings are either even or odd and are numbered  by ascending sequences of even or odd numbers, as they {\bf first} show up as one travels along the knot from the starting point. The Gauss code is then the sequence of the numbers of the crossings as they show up from the starting point onwards with a positive sign if one passes over the crossing point and a negative sign otherwise. In particular, a Gauss code is a list of pairs of integers of the same absolute value and oppposite signs (see \cite{CetS}).

Consider a Gauss code of a planar diagram of a particular knot. Successive negative integers in the Gauss code stand for the arcs of the planar diagram and we can think of the positive numbers as standing for the crossings. Pick one crossing. The coloring condition to be read off here is of the form $a\triangleright b=c$. The parity (even or odd) of the corresponding positive integer (which lies between successive negative integers) tells one right away which arcs are the first factor of the left hand-side and the right hand-side of the coloring condition associated to this crossing. Repeating this procedure for each positive integer, i.e., for each crossing, we set up a system of equations whose solutions are the colorings for the knot.

The methods referred to above were implemented with C and the results agreed.
The implementation using Burau matrices (eigenvector equation) consists mainly in selecting from all candidate colorings those that satisfy the equation. As for the implementation based on Gauss codes, since it implies a much larger number of candidate colorings, a prior selection of relevant arcs was done and the candidate colorings were restricted to those arcs.

The matrix with the invariants is displayed in Appendix A. The list of knots used in this paper is ordered as found in \cite{J}. The list of quandles used is $R_3$, $R_5$, $R_7$, ${\mathbb Z}_7[T, T^{-1}]/(T-2)$, ${\mathbb Z}_3[T, T^{-1}]/(T^2+1)$, ${\mathbb Z}_2[T, T^{-1}]/(T^2+T+1)$, ${\mathbb Z}_3[T, T^{-1}]/(T^2+T-1)$, ${\mathbb Z}_2[T, T^{-1}]/(T^3+T^2+1)$, ${\mathbb Z}_5[T, T^{-1}]/(T-2)$, $R_9$ and was set up after some experiments to see which quandles were the most efficient ones. In this way, the entry ``6'' indexed by ``17'' on the first column and ``48'' on the first row, tells one that the quandle ${\mathbb Z}_2[T, T^{-1}]/(T^2+T+1)$ is the first one from the list of quandles that colors knots $8_3$ and $9_{12}$ in different numbers of ways whereas entry ``0'' indexed by ``90'' on the first column and ``95'' on the first row, tells one that knots $10_6$ and $10_{11}$ are ${\bf not}$ colored in different numbers of ways by any of the quandles in the list above.

There were 793 inconclusive cases (zeros) in the 30876 cases studied. Using Fenn and Rourke's criterion (only trivial colorings for one knot and at least one non-trivial coloring for the other one, see \cite{FR}) there were 962 inconclusive cases.

In Fenn and Rourke's \cite{FR} it is conjectured that there should exist a countable sequence of finite racks such that the sequence of sets of homomorphisms from the fundamental rack of a given knot to the racks in the sequence is a complete invariant of the knot. Since quandles are specializations of racks and despite the fact that we have merely touched on an insignificant part of all (classical) knots, we feel we have given strength to that conjecture.

%Ncomparacoes: 30876 Nzeros: 793 Nzeros triv:962

\clearpage

\appendix
\section{The Matrix}
\vskip 100pt
\begin{table}[h]\tiny{
\begin{center}
\begin{tabular}{p{.35cm}|*{20}{p{.35cm}}}
    & 1 & 2 & 3 & 4 & 5 & 6 & 7 & 8 & 9 & 10 & 11 & 12 & 13 & 14 & 15 & 16 & 17 & 18 & 19 & 20 \\\hline 
1   &   & 1 & 1 & 1 & 4 & 1 & 1 & 1 & 1 & 1 & 2 & 1 & 1 & 3 & 1 & 1 & 1 & 1 & 3 & 1 \\
2   &   &   & 5 & 2 & 1 & 2 & 2 & 2 & 2 & 2 & 1 & 2 & 2 & 1 & 2 & 2 & 2 & 2 & 1 & 2 \\
3   &   &   &   & 2 & 1 & 2 & 2 & 2 & 2 & 2 & 1 & 2 & 2 & 1 & 2 & 2 & 2 & 2 & 1 & 2 \\
4   &   &   &   &   & 1 & 3 & 3 & 5 & 3 & 3 & 1 & 3 & 3 & 1 & 3 & 3 & 3 & 3 & 1 & 3 \\
5   &   &   &   &   &   & 1 & 1 & 1 & 1 & 1 & 2 & 1 & 1 & 3 & 1 & 1 & 1 & 1 & 3 & 1 \\
6   &   &   &   &   &   &   & 4 & 3 & 6 & 6 & 1 & 4 & 7 & 1 & 6 & 5 & 5 & 6 & 1 & 5 \\
7   &   &   &   &   &   &   &   & 3 & 4 & 4 & 1 & 5 & 4 & 1 & 4 & 4 & 4 & 4 & 1 & 4 \\
8   &   &   &   &   &   &   &   &   & 3 & 3 & 1 & 3 & 3 & 1 & 3 & 3 & 3 & 3 & 1 & 3 \\
9   &   &   &   &   &   &   &   &   &   & 7 & 1 & 4 & 6 & 1 & 10 & 5 & 5 & 10 & 1 & 5 \\
10   &   &   &   &   &   &   &   &   &   &   & 1 & 4 & 6 & 1 & 7 & 5 & 5 & 7 & 1 & 5 \\
11   &   &   &   &   &   &   &   &   &   &   &   & 1 & 1 & 2 & 1 & 1 & 1 & 1 & 2 & 1 \\
12   &   &   &   &   &   &   &   &   &   &   &   &   & 4 & 1 & 4 & 4 & 4 & 4 & 1 & 4 \\
13   &   &   &   &   &   &   &   &   &   &   &   &   &   & 1 & 6 & 5 & 5 & 6 & 1 & 5 \\
14   &   &   &   &   &   &   &   &   &   &   &   &   &   &   & 1 & 1 & 1 & 1 & 6 & 1 \\
15   &   &   &   &   &   &   &   &   &   &   &   &   &   &   &   & 5 & 5 & 0 & 1 & 5 \\
16   &   &   &   &   &   &   &   &   &   &   &   &   &   &   &   &   & 8 & 5 & 1 & 8 \\
17   &   &   &   &   &   &   &   &   &   &   &   &   &   &   &   &   &   & 5 & 1 & 0 \\
18   &   &   &   &   &   &   &   &   &   &   &   &   &   &   &   &   &   &   & 1 & 5 \\
19   &   &   &   &   &   &   &   &   &   &   &   &   &   &   &   &   &   &   &   & 1 \\
20   &   &   &   &   &   &   &   &   &   &   &   &   &   &   &   &   &   &   &   &   
\end{tabular}
\vskip 20pt
{\large Table 1}
\end{center}
 }
\end{table}

\begin{table}[h]\tiny{
\begin{center}
\begin{tabular}{p{.35cm}|*{20}{p{.35cm}}}
    & 21 & 22 & 23 & 24 & 25 & 26 & 27 & 28 & 29 & 30 & 31 & 32 & 33 & 34 & 35 & 36 & 37 & 38 & 39 & 40 \\\hline 
1   & 1 & 1 & 1 & 9 & 4 & 1 & 1 & 1 & 10 & 1 & 1 & 1 & 5 & 9 & 2 & 9 & 2 & 1 & 3 & 1 \\
2   & 2 & 5 & 4 & 1 & 1 & 2 & 2 & 2 & 1 & 3 & 2 & 1 & 1 & 1 & 1 & 1 & 1 & 2 & 1 & 2 \\
3   & 2 & 7 & 4 & 1 & 1 & 2 & 2 & 2 & 1 & 3 & 2 & 1 & 1 & 1 & 1 & 1 & 1 & 2 & 1 & 2 \\
4   & 3 & 2 & 2 & 1 & 1 & 3 & 3 & 3 & 1 & 2 & 3 & 1 & 1 & 1 & 1 & 1 & 1 & 3 & 1 & 3 \\
5   & 1 & 1 & 1 & 4 & 6 & 1 & 1 & 1 & 4 & 1 & 1 & 1 & 4 & 4 & 2 & 4 & 2 & 1 & 3 & 1 \\
6   & 7 & 2 & 2 & 1 & 1 & 7 & 6 & 7 & 1 & 2 & 5 & 1 & 1 & 1 & 1 & 1 & 1 & 5 & 1 & 7 \\
7   & 4 & 2 & 2 & 1 & 1 & 4 & 4 & 4 & 1 & 2 & 4 & 1 & 1 & 1 & 1 & 1 & 1 & 4 & 1 & 4 \\
\end{tabular}
\vskip 20pt
{\large Table 2}
\end{center}
 }
\end{table}\clearpage

\begin{table}\tiny{
\begin{center}
\begin{tabular}{p{.35cm}|*{19}{p{.35cm}}}
    & 41 & 42 & 43 & 44 & 45 & 46 & 47 & 48 & 49 & 50 & 51 & 52 & 53 & 54 & 55 & 56 & 57 & 58 & 59 \\\hline 
1   & 5 & 1 & 1 & 1 & 6 & 4 & 1 & 1 & 1 & 6 & 7 & 6 & 1 & 1 & 1 & 1 & 1 & 2 & 2 \\
2   & 1 & 2 & 2 & 2 & 1 & 1 & 3 & 2 & 2 & 1 & 1 & 1 & 2 & 2 & 2 & 2 & 2 & 1 & 1 \\
3   & 1 & 2 & 2 & 2 & 1 & 1 & 3 & 2 & 2 & 1 & 1 & 1 & 2 & 2 & 2 & 2 & 2 & 1 & 1 \\
4   & 1 & 3 & 3 & 3 & 1 & 1 & 2 & 3 & 3 & 1 & 1 & 1 & 3 & 3 & 3 & 3 & 3 & 1 & 1 \\
5   & 4 & 1 & 1 & 1 & 4 & 5 & 1 & 1 & 1 & 4 & 4 & 4 & 1 & 1 & 1 & 1 & 1 & 2 & 2 \\
6   & 1 & 5 & 7 & 7 & 1 & 1 & 2 & 5 & 6 & 1 & 1 & 1 & 7 & 7 & 7 & 6 & 4 & 1 & 1 \\
7   & 1 & 4 & 4 & 4 & 1 & 1 & 2 & 4 & 4 & 1 & 1 & 1 & 4 & 4 & 4 & 4 & 5 & 1 & 1 \\
\end{tabular}
\vskip 20pt
{\large Table 3}
\end{center}
 }
\end{table}\clearpage

\begin{table}\tiny{
\begin{center}
\begin{tabular}{p{.35cm}|*{19}{p{.35cm}}}
    & 60 & 61 & 62 & 63 & 64 & 65 & 66 & 67 & 68 & 69 & 70 & 71 & 72 & 73 & 74 & 75 & 76 & 77 & 78 \\\hline 
1   & 1 & 1 & 1 & 10 & 10 & 1 & 1 & 1 & 1 & 6 & 1 & 1 & 1 & 4 & 1 & 2 & 1 & 1 & 1 \\
2   & 2 & 2 & 2 & 1 & 1 & 2 & 5 & 2 & 2 & 1 & 1 & 2 & 1 & 1 & 10 & 1 & 2 & 2 & 2 \\
3   & 2 & 2 & 2 & 1 & 1 & 2 & 8 & 2 & 2 & 1 & 1 & 2 & 1 & 1 & 5 & 1 & 2 & 2 & 2 \\
4   & 3 & 3 & 5 & 1 & 1 & 3 & 2 & 3 & 3 & 1 & 1 & 3 & 1 & 1 & 2 & 1 & 3 & 5 & 3 \\
5   & 1 & 1 & 1 & 4 & 4 & 1 & 1 & 1 & 1 & 4 & 1 & 1 & 1 & 5 & 1 & 2 & 1 & 1 & 1 \\
6   & 6 & 5 & 3 & 1 & 1 & 6 & 2 & 7 & 7 & 1 & 1 & 6 & 1 & 1 & 2 & 1 & 3 & 3 & 6 \\
7   & 4 & 4 & 3 & 1 & 1 & 4 & 2 & 4 & 4 & 1 & 1 & 4 & 1 & 1 & 2 & 1 & 3 & 3 & 4 \\
\end{tabular}
\vskip 20pt
{\large Table 4}
\end{center}
 }
\end{table}\clearpage

\begin{table}\tiny{
\begin{center}
\begin{tabular}{p{.35cm}|*{19}{p{.35cm}}}
    & 79 & 80 & 81 & 82 & 83 & 84 & 85 & 86 & 87 & 88 & 89 & 90 & 91 & 92 & 93 & 94 & 95 & 96 & 97 \\\hline 
1   & 1 & 1 & 1 & 1 & 1 & 1 & 1 & 1 & 1 & 6 & 5 & 1 & 1 & 1 & 0 & 2 & 1 & 1 & 1 \\
2   & 2 & 2 & 1 & 1 & 1 & 2 & 2 & 2 & 5 & 1 & 1 & 2 & 2 & 2 & 1 & 1 & 2 & 2 & 2 \\
3   & 2 & 2 & 1 & 1 & 1 & 2 & 2 & 2 & 0 & 1 & 1 & 2 & 2 & 2 & 1 & 1 & 2 & 2 & 2 \\
4   & 3 & 3 & 1 & 1 & 1 & 2 & 3 & 3 & 2 & 1 & 1 & 3 & 3 & 3 & 1 & 1 & 3 & 3 & 3 \\
5   & 1 & 1 & 1 & 1 & 1 & 1 & 1 & 1 & 1 & 4 & 4 & 1 & 1 & 1 & 4 & 2 & 1 & 1 & 1 \\
6   & 6 & 4 & 1 & 1 & 1 & 2 & 5 & 6 & 2 & 1 & 1 & 6 & 5 & 7 & 1 & 1 & 6 & 7 & 6 \\
7   & 4 & 5 & 1 & 1 & 1 & 2 & 4 & 4 & 2 & 1 & 1 & 4 & 4 & 4 & 1 & 1 & 4 & 4 & 4 \\
\end{tabular}
\vskip 20pt
{\large Table 5}
\end{center}
 }
\end{table}\clearpage

\begin{table}\tiny{
\begin{center}
\begin{tabular}{p{.35cm}|*{19}{p{.35cm}}}
    & 79 & 80 & 81 & 82 & 83 & 84 & 85 & 86 & 87 & 88 & 89 & 90 & 91 & 92 & 93 & 94 & 95 & 96 & 97 \\\hline 
84   &   &   &   &   &   &   & 2 & 2 & 2 & 1 & 1 & 2 & 2 & 2 & 1 & 1 & 2 & 2 & 2 \\
85   &   &   &   &   &   &   &   & 5 & 2 & 1 & 1 & 5 & 0 & 5 & 1 & 1 & 5 & 5 & 5 \\
86   &   &   &   &   &   &   &   &   & 2 & 1 & 1 & 10 & 5 & 6 & 1 & 1 & 10 & 6 & 0 \\
87   &   &   &   &   &   &   &   &   &   & 1 & 1 & 2 & 2 & 2 & 1 & 1 & 2 & 2 & 2 \\
88   &   &   &   &   &   &   &   &   &   &   & 5 & 1 & 1 & 1 & 6 & 2 & 1 & 1 & 1 \\
89   &   &   &   &   &   &   &   &   &   &   &   & 1 & 1 & 1 & 5 & 2 & 1 & 1 & 1 \\
90   &   &   &   &   &   &   &   &   &   &   &   &   & 5 & 6 & 1 & 1 & 0 & 6 & 10 \\
91   &   &   &   &   &   &   &   &   &   &   &   &   &   & 5 & 1 & 1 & 5 & 5 & 5 \\
92   &   &   &   &   &   &   &   &   &   &   &   &   &   &   & 1 & 1 & 6 & 10 & 6 \\
93   &   &   &   &   &   &   &   &   &   &   &   &   &   &   &   & 2 & 1 & 1 & 1 \\
94   &   &   &   &   &   &   &   &   &   &   &   &   &   &   &   &   & 1 & 1 & 1 \\
95   &   &   &   &   &   &   &   &   &   &   &   &   &   &   &   &   &   & 6 & 10 \\
96   &   &   &   &   &   &   &   &   &   &   &   &   &   &   &   &   &   &   & 6 \\
97   &   &   &   &   &   &   &   &   &   &   &   &   &   &   &   &   &   &   &   
\end{tabular}
\vskip 20pt
{\large Table 6}
\end{center}
 }
\end{table}\clearpage

\begin{table}\tiny{
\begin{center}
\begin{tabular}{p{.35cm}|*{19}{p{.35cm}}}
    & 98 & 99 & 100 & 101 & 102 & 103 & 104 & 105 & 106 & 107 & 108 & 109 & 110 & 111 & 112 & 113 & 114 & 115 & 116 \\\hline 
1   & 5 & 1 & 1 & 1 & 1 & 5 & 1 & 2 & 1 & 1 & 1 & 1 & 1 & 1 & 1 & 3 & 1 & 6 & 5 \\
2   & 1 & 2 & 2 & 2 & 5 & 1 & 3 & 1 & 2 & 2 & 5 & 5 & 2 & 2 & 2 & 1 & 2 & 1 & 1 \\
3   & 1 & 2 & 2 & 2 & 0 & 1 & 3 & 1 & 2 & 2 & 0 & 0 & 2 & 2 & 2 & 1 & 2 & 1 & 1 \\
4   & 1 & 3 & 3 & 3 & 2 & 1 & 2 & 1 & 5 & 3 & 2 & 2 & 3 & 3 & 3 & 1 & 3 & 1 & 1 \\
5   & 4 & 1 & 1 & 1 & 1 & 4 & 1 & 2 & 1 & 1 & 1 & 1 & 1 & 1 & 1 & 3 & 1 & 4 & 4 \\
6   & 1 & 5 & 0 & 6 & 2 & 1 & 2 & 1 & 3 & 7 & 2 & 2 & 5 & 5 & 4 & 1 & 6 & 1 & 1 \\
7   & 1 & 4 & 4 & 4 & 2 & 1 & 2 & 1 & 3 & 4 & 2 & 2 & 4 & 4 & 5 & 1 & 4 & 1 & 1 \\
\end{tabular}
\vskip 20pt
{\large Table 7}
\end{center}
 }
\end{table}\clearpage

\begin{table}\tiny{
\begin{center}
\begin{tabular}{p{.35cm}|*{19}{p{.35cm}}}
    & 117 & 118 & 119 & 120 & 121 & 122 & 123 & 124 & 125 & 126 & 127 & 128 & 129 & 130 & 131 & 132 & 133 & 134 & 135 \\\hline 
1   & 1 & 1 & 1 & 6 & 1 & 1 & 1 & 2 & 1 & 6 & 1 & 1 & 1 & 1 & 1 & 1 & 1 & 1 & 1 \\
2   & 5 & 2 & 2 & 1 & 2 & 2 & 2 & 1 & 2 & 1 & 2 & 2 & 2 & 2 & 2 & 2 & 2 & 2 & 2 \\
3   & 0 & 2 & 2 & 1 & 2 & 2 & 2 & 1 & 2 & 1 & 2 & 2 & 2 & 2 & 2 & 2 & 2 & 2 & 2 \\
4   & 2 & 3 & 5 & 1 & 3 & 3 & 3 & 1 & 3 & 1 & 3 & 3 & 3 & 3 & 3 & 4 & 3 & 3 & 3 \\
5   & 1 & 1 & 1 & 4 & 1 & 1 & 1 & 2 & 1 & 4 & 1 & 1 & 1 & 1 & 1 & 1 & 1 & 1 & 1 \\
6   & 2 & 7 & 3 & 1 & 4 & 6 & 6 & 1 & 7 & 1 & 5 & 4 & 7 & 7 & 5 & 3 & 5 & 7 & 7 \\
7   & 2 & 4 & 3 & 1 & 5 & 4 & 4 & 1 & 4 & 1 & 4 & 5 & 4 & 4 & 4 & 3 & 4 & 4 & 4 \\
\end{tabular}
\vskip 20pt
{\large Table 8}
\end{center}
 }
\end{table}\clearpage

\begin{table}\tiny{
\begin{center}
\begin{tabular}{p{.35cm}|*{19}{p{.35cm}}}
    & 117 & 118 & 119 & 120 & 121 & 122 & 123 & 124 & 125 & 126 & 127 & 128 & 129 & 130 & 131 & 132 & 133 & 134 & 135 \\\hline 
84   & 2 & 2 & 2 & 1 & 2 & 2 & 2 & 1 & 2 & 1 & 2 & 2 & 2 & 2 & 2 & 2 & 2 & 2 & 2 \\
85   & 2 & 5 & 3 & 1 & 4 & 5 & 5 & 1 & 5 & 1 & 8 & 4 & 5 & 5 & 0 & 3 & 0 & 5 & 5 \\
86   & 2 & 6 & 3 & 1 & 4 & 7 & 0 & 1 & 6 & 1 & 5 & 4 & 6 & 6 & 5 & 3 & 5 & 6 & 6 \\
87   & 0 & 2 & 2 & 1 & 2 & 2 & 2 & 1 & 2 & 1 & 2 & 2 & 2 & 2 & 2 & 2 & 2 & 2 & 2 \\
88   & 1 & 1 & 1 & 9 & 1 & 1 & 1 & 2 & 1 & 7 & 1 & 1 & 1 & 1 & 1 & 1 & 1 & 1 & 1 \\
89   & 1 & 1 & 1 & 5 & 1 & 1 & 1 & 2 & 1 & 5 & 1 & 1 & 1 & 1 & 1 & 1 & 1 & 1 & 1 \\
90   & 2 & 6 & 3 & 1 & 4 & 7 & 10 & 1 & 6 & 1 & 5 & 4 & 6 & 6 & 5 & 3 & 5 & 6 & 6 \\
\end{tabular}
\vskip 20pt
{\large Table 9}
\end{center}
 }
\end{table}\clearpage

\begin{table}\tiny{
\begin{center}
\begin{tabular}{p{.35cm}|*{19}{p{.35cm}}}
    & 136 & 137 & 138 & 139 & 140 & 141 & 142 & 143 & 144 & 145 & 146 & 147 & 148 & 149 & 150 & 151 & 152 & 153 & 154 \\\hline 
1   & 1 & 1 & 1 & 1 & 1 & 1 & 1 & 2 & 1 & 5 & 2 & 4 & 0 & 3 & 2 & 3 & 6 & 5 & 1 \\
2   & 2 & 2 & 2 & 2 & 5 & 2 & 0 & 1 & 10 & 1 & 1 & 1 & 1 & 1 & 1 & 1 & 1 & 1 & 2 \\
3   & 2 & 2 & 2 & 2 & 0 & 2 & 5 & 1 & 5 & 1 & 1 & 1 & 1 & 1 & 1 & 1 & 1 & 1 & 2 \\
4   & 3 & 3 & 3 & 3 & 2 & 3 & 2 & 1 & 2 & 1 & 1 & 1 & 1 & 1 & 1 & 1 & 1 & 1 & 3 \\
5   & 1 & 1 & 1 & 1 & 1 & 1 & 1 & 2 & 1 & 4 & 2 & 5 & 4 & 3 & 2 & 3 & 4 & 4 & 1 \\
6   & 7 & 7 & 7 & 0 & 2 & 7 & 2 & 1 & 2 & 1 & 1 & 1 & 1 & 1 & 1 & 1 & 1 & 1 & 6 \\
7   & 4 & 4 & 4 & 4 & 2 & 4 & 2 & 1 & 2 & 1 & 1 & 1 & 1 & 1 & 1 & 1 & 1 & 1 & 4 \\
\end{tabular}
\vskip 20pt
{\large Table 10}
\end{center}
 }
\end{table}\clearpage

\begin{table}\tiny{
\begin{center}
\begin{tabular}{p{.35cm}|*{19}{p{.35cm}}}
    & 136 & 137 & 138 & 139 & 140 & 141 & 142 & 143 & 144 & 145 & 146 & 147 & 148 & 149 & 150 & 151 & 152 & 153 & 154 \\\hline 
84   & 2 & 2 & 2 & 2 & 2 & 2 & 2 & 1 & 2 & 1 & 1 & 1 & 1 & 1 & 1 & 1 & 1 & 1 & 2 \\
85   & 5 & 5 & 5 & 5 & 2 & 5 & 2 & 1 & 2 & 1 & 1 & 1 & 1 & 1 & 1 & 1 & 1 & 1 & 5 \\
86   & 6 & 6 & 6 & 6 & 2 & 6 & 2 & 1 & 2 & 1 & 1 & 1 & 1 & 1 & 1 & 1 & 1 & 1 & 10 \\
87   & 2 & 2 & 2 & 2 & 0 & 2 & 5 & 1 & 5 & 1 & 1 & 1 & 1 & 1 & 1 & 1 & 1 & 1 & 2 \\
88   & 1 & 1 & 1 & 1 & 1 & 1 & 1 & 2 & 1 & 5 & 2 & 4 & 6 & 3 & 2 & 3 & 9 & 5 & 1 \\
89   & 1 & 1 & 1 & 1 & 1 & 1 & 1 & 2 & 1 & 6 & 2 & 4 & 5 & 3 & 2 & 3 & 5 & 6 & 1 \\
90   & 6 & 6 & 6 & 6 & 2 & 6 & 2 & 1 & 2 & 1 & 1 & 1 & 1 & 1 & 1 & 1 & 1 & 1 & 0 \\
\end{tabular}
\vskip 20pt
{\large Table 11}
\end{center}
 }
\end{table}\clearpage

\begin{table}\tiny{
\begin{center}
\begin{tabular}{p{.35cm}|*{19}{p{.35cm}}}
    & 155 & 156 & 157 & 158 & 159 & 160 & 161 & 162 & 163 & 164 & 165 & 166 & 167 & 168 & 169 & 170 & 171 & 172 & 173 \\\hline 
1   & 1 & 1 & 1 & 1 & 1 & 0 & 3 & 4 & 1 & 1 & 1 & 3 & 1 & 0 & 5 & 1 & 4 & 1 & 6 \\
2   & 2 & 2 & 2 & 1 & 1 & 1 & 1 & 1 & 2 & 2 & 5 & 1 & 5 & 1 & 1 & 2 & 1 & 2 & 1 \\
3   & 2 & 2 & 2 & 1 & 1 & 1 & 1 & 1 & 2 & 2 & 0 & 1 & 0 & 1 & 1 & 2 & 1 & 2 & 1 \\
4   & 5 & 3 & 3 & 1 & 1 & 1 & 1 & 1 & 3 & 3 & 2 & 1 & 2 & 1 & 1 & 3 & 1 & 3 & 1 \\
5   & 1 & 1 & 1 & 1 & 1 & 4 & 3 & 6 & 1 & 1 & 1 & 3 & 1 & 4 & 4 & 1 & 6 & 1 & 4 \\
6   & 3 & 6 & 6 & 1 & 1 & 1 & 1 & 1 & 5 & 7 & 2 & 1 & 2 & 1 & 1 & 7 & 1 & 4 & 1 \\
7   & 3 & 4 & 4 & 1 & 1 & 1 & 1 & 1 & 4 & 4 & 2 & 1 & 2 & 1 & 1 & 4 & 1 & 5 & 1 \\
\end{tabular}
\vskip 20pt
{\large Table 12}
\end{center}
 }
\end{table}\clearpage

\begin{table}\tiny{
\begin{center}
\begin{tabular}{p{.35cm}|*{19}{p{.35cm}}}
    & 155 & 156 & 157 & 158 & 159 & 160 & 161 & 162 & 163 & 164 & 165 & 166 & 167 & 168 & 169 & 170 & 171 & 172 & 173 \\\hline 
84   & 2 & 2 & 2 & 1 & 1 & 1 & 1 & 1 & 2 & 2 & 2 & 1 & 2 & 1 & 1 & 2 & 1 & 2 & 1 \\
85   & 3 & 5 & 5 & 1 & 1 & 1 & 1 & 1 & 0 & 5 & 2 & 1 & 2 & 1 & 1 & 5 & 1 & 4 & 1 \\
86   & 3 & 0 & 0 & 1 & 1 & 1 & 1 & 1 & 5 & 6 & 2 & 1 & 2 & 1 & 1 & 6 & 1 & 4 & 1 \\
87   & 2 & 2 & 2 & 1 & 1 & 1 & 1 & 1 & 2 & 2 & 0 & 1 & 0 & 1 & 1 & 2 & 1 & 2 & 1 \\
88   & 1 & 1 & 1 & 1 & 1 & 6 & 3 & 4 & 1 & 1 & 1 & 3 & 1 & 6 & 5 & 1 & 4 & 1 & 0 \\
89   & 1 & 1 & 1 & 1 & 1 & 5 & 3 & 4 & 1 & 1 & 1 & 3 & 1 & 5 & 7 & 1 & 4 & 1 & 5 \\
90   & 3 & 10 & 10 & 1 & 1 & 1 & 1 & 1 & 5 & 6 & 2 & 1 & 2 & 1 & 1 & 6 & 1 & 4 & 1 \\
\end{tabular}
\vskip 20pt
{\large Table 13}
\end{center}
 }
\end{table}\clearpage

\begin{table}\tiny{
\begin{center}
\begin{tabular}{p{.35cm}|*{19}{p{.35cm}}}
    & 155 & 156 & 157 & 158 & 159 & 160 & 161 & 162 & 163 & 164 & 165 & 166 & 167 & 168 & 169 & 170 & 171 & 172 & 173 \\\hline 
167   &   &   &   &   &   &   &   &   &   &   &   &   &   & 1 & 1 & 2 & 1 & 2 & 1 \\
168   &   &   &   &   &   &   &   &   &   &   &   &   &   &   & 5 & 1 & 4 & 1 & 6 \\
169   &   &   &   &   &   &   &   &   &   &   &   &   &   &   &   & 1 & 4 & 1 & 5 \\
170   &   &   &   &   &   &   &   &   &   &   &   &   &   &   &   &   & 1 & 4 & 1 \\
171   &   &   &   &   &   &   &   &   &   &   &   &   &   &   &   &   &   & 1 & 4 \\
172   &   &   &   &   &   &   &   &   &   &   &   &   &   &   &   &   &   &   & 1 \\
173   &   &   &   &   &   &   &   &   &   &   &   &   &   &   &   &   &   &   &   
\end{tabular}
\vskip 20pt
{\large Table 14}
\end{center}
 }
\end{table}\clearpage

\begin{table}\tiny{
\begin{center}
\begin{tabular}{p{.35cm}|*{19}{p{.35cm}}}
    & 174 & 175 & 176 & 177 & 178 & 179 & 180 & 181 & 182 & 183 & 184 & 185 & 186 & 187 & 188 & 189 & 190 & 191 & 192 \\\hline 
1   & 1 & 1 & 1 & 1 & 1 & 1 & 7 & 0 & 1 & 1 & 1 & 1 & 1 & 2 & 1 & 1 & 2 & 5 & 3 \\
2   & 2 & 2 & 2 & 2 & 2 & 2 & 1 & 1 & 1 & 1 & 4 & 6 & 2 & 1 & 2 & 2 & 1 & 1 & 1 \\
3   & 2 & 2 & 2 & 2 & 2 & 2 & 1 & 1 & 1 & 1 & 4 & 5 & 2 & 1 & 2 & 2 & 1 & 1 & 1 \\
4   & 4 & 3 & 3 & 3 & 3 & 4 & 1 & 1 & 1 & 1 & 2 & 2 & 3 & 1 & 6 & 5 & 1 & 1 & 1 \\
5   & 1 & 1 & 1 & 1 & 1 & 1 & 4 & 4 & 1 & 1 & 1 & 1 & 1 & 2 & 1 & 1 & 2 & 4 & 3 \\
6   & 3 & 7 & 7 & 6 & 7 & 3 & 1 & 1 & 1 & 1 & 2 & 2 & 7 & 1 & 3 & 3 & 1 & 1 & 1 \\
7   & 3 & 4 & 4 & 4 & 4 & 3 & 1 & 1 & 1 & 1 & 2 & 2 & 4 & 1 & 3 & 3 & 1 & 1 & 1 \\
\end{tabular}
\vskip 20pt
{\large Table 15}
\end{center}
 }
\end{table}\clearpage

\begin{table}\tiny{
\begin{center}
\begin{tabular}{p{.35cm}|*{19}{p{.35cm}}}
    & 174 & 175 & 176 & 177 & 178 & 179 & 180 & 181 & 182 & 183 & 184 & 185 & 186 & 187 & 188 & 189 & 190 & 191 & 192 \\\hline 
84   & 2 & 2 & 2 & 2 & 2 & 2 & 1 & 1 & 1 & 1 & 2 & 2 & 2 & 1 & 2 & 2 & 1 & 1 & 1 \\
85   & 3 & 5 & 5 & 5 & 5 & 3 & 1 & 1 & 1 & 1 & 2 & 2 & 5 & 1 & 3 & 3 & 1 & 1 & 1 \\
86   & 3 & 6 & 6 & 0 & 6 & 3 & 1 & 1 & 1 & 1 & 2 & 2 & 6 & 1 & 3 & 3 & 1 & 1 & 1 \\
87   & 2 & 2 & 2 & 2 & 2 & 2 & 1 & 1 & 1 & 1 & 4 & 5 & 2 & 1 & 2 & 2 & 1 & 1 & 1 \\
88   & 1 & 1 & 1 & 1 & 1 & 1 & 6 & 6 & 1 & 1 & 1 & 1 & 1 & 2 & 1 & 1 & 2 & 5 & 3 \\
89   & 1 & 1 & 1 & 1 & 1 & 1 & 5 & 5 & 1 & 1 & 1 & 1 & 1 & 2 & 1 & 1 & 2 & 6 & 3 \\
90   & 3 & 6 & 6 & 10 & 6 & 3 & 1 & 1 & 1 & 1 & 2 & 2 & 6 & 1 & 3 & 3 & 1 & 1 & 1 \\
\end{tabular}
\vskip 20pt
{\large Table 16}
\end{center}
 }
\end{table}\clearpage

\begin{table}\tiny{
\begin{center}
\begin{tabular}{p{.35cm}|*{19}{p{.35cm}}}
    & 174 & 175 & 176 & 177 & 178 & 179 & 180 & 181 & 182 & 183 & 184 & 185 & 186 & 187 & 188 & 189 & 190 & 191 & 192 \\\hline 
167   & 2 & 2 & 2 & 2 & 2 & 2 & 1 & 1 & 1 & 1 & 4 & 5 & 2 & 1 & 2 & 2 & 1 & 1 & 1 \\
168   & 1 & 1 & 1 & 1 & 1 & 1 & 7 & 0 & 1 & 1 & 1 & 1 & 1 & 2 & 1 & 1 & 2 & 5 & 3 \\
169   & 1 & 1 & 1 & 1 & 1 & 1 & 5 & 5 & 1 & 1 & 1 & 1 & 1 & 2 & 1 & 1 & 2 & 6 & 3 \\
170   & 3 & 10 & 0 & 6 & 10 & 3 & 1 & 1 & 1 & 1 & 2 & 2 & 0 & 1 & 3 & 3 & 1 & 1 & 1 \\
171   & 1 & 1 & 1 & 1 & 1 & 1 & 4 & 4 & 1 & 1 & 1 & 1 & 1 & 2 & 1 & 1 & 2 & 4 & 3 \\
172   & 3 & 4 & 4 & 4 & 4 & 3 & 1 & 1 & 1 & 1 & 2 & 2 & 4 & 1 & 3 & 3 & 1 & 1 & 1 \\
173   & 1 & 1 & 1 & 1 & 1 & 1 & 6 & 6 & 1 & 1 & 1 & 1 & 1 & 2 & 1 & 1 & 2 & 5 & 3 \\
\end{tabular}
\vskip 20pt
{\large Table 17}
\end{center}
 }
\end{table}\clearpage

\begin{table}\tiny{
\begin{center}
\begin{tabular}{p{.35cm}|*{19}{p{.35cm}}}
    & 193 & 194 & 195 & 196 & 197 & 198 & 199 & 200 & 201 & 202 & 203 & 204 & 205 & 206 & 207 & 208 & 209 & 210 & 211 \\\hline 
1   & 1 & 1 & 1 & 4 & 7 & 0 & 1 & 1 & 1 & 1 & 1 & 2 & 1 & 2 & 1 & 1 & 1 & 1 & 1 \\
2   & 6 & 2 & 2 & 1 & 1 & 1 & 2 & 5 & 2 & 2 & 2 & 1 & 5 & 1 & 2 & 2 & 2 & 2 & 2 \\
3   & 5 & 2 & 2 & 1 & 1 & 1 & 2 & 0 & 2 & 2 & 2 & 1 & 0 & 1 & 2 & 2 & 2 & 2 & 2 \\
4   & 2 & 3 & 7 & 1 & 1 & 1 & 3 & 2 & 3 & 3 & 3 & 1 & 2 & 1 & 3 & 3 & 3 & 3 & 3 \\
5   & 1 & 1 & 1 & 6 & 4 & 4 & 1 & 1 & 1 & 1 & 1 & 2 & 1 & 2 & 1 & 1 & 1 & 1 & 1 \\
6   & 2 & 4 & 3 & 1 & 1 & 1 & 6 & 2 & 7 & 7 & 6 & 1 & 2 & 1 & 7 & 7 & 5 & 7 & 4 \\
7   & 2 & 0 & 3 & 1 & 1 & 1 & 4 & 2 & 4 & 4 & 4 & 1 & 2 & 1 & 4 & 4 & 4 & 4 & 5 \\
\end{tabular}
\vskip 20pt
{\large Table 18}
\end{center}
 }
\end{table}\clearpage

\begin{table}\tiny{
\begin{center}
\begin{tabular}{p{.35cm}|*{19}{p{.35cm}}}
    & 193 & 194 & 195 & 196 & 197 & 198 & 199 & 200 & 201 & 202 & 203 & 204 & 205 & 206 & 207 & 208 & 209 & 210 & 211 \\\hline 
84   & 2 & 2 & 2 & 1 & 1 & 1 & 2 & 2 & 2 & 2 & 2 & 1 & 2 & 1 & 2 & 2 & 2 & 2 & 2 \\
85   & 2 & 4 & 3 & 1 & 1 & 1 & 5 & 2 & 5 & 5 & 5 & 1 & 2 & 1 & 5 & 5 & 0 & 5 & 4 \\
86   & 2 & 4 & 3 & 1 & 1 & 1 & 6 & 2 & 6 & 6 & 0 & 1 & 2 & 1 & 6 & 6 & 5 & 6 & 4 \\
87   & 5 & 2 & 2 & 1 & 1 & 1 & 2 & 0 & 2 & 2 & 2 & 1 & 0 & 1 & 2 & 2 & 2 & 2 & 2 \\
88   & 1 & 1 & 1 & 4 & 6 & 6 & 1 & 1 & 1 & 1 & 1 & 2 & 1 & 2 & 1 & 1 & 1 & 1 & 1 \\
89   & 1 & 1 & 1 & 4 & 5 & 5 & 1 & 1 & 1 & 1 & 1 & 2 & 1 & 2 & 1 & 1 & 1 & 1 & 1 \\
90   & 2 & 4 & 3 & 1 & 1 & 1 & 6 & 2 & 6 & 6 & 10 & 1 & 2 & 1 & 6 & 6 & 5 & 6 & 4 \\
\end{tabular}
\vskip 20pt
{\large Table 19}
\end{center}
 }
\end{table}\clearpage

\begin{table}\tiny{
\begin{center}
\begin{tabular}{p{.35cm}|*{19}{p{.35cm}}}
    & 193 & 194 & 195 & 196 & 197 & 198 & 199 & 200 & 201 & 202 & 203 & 204 & 205 & 206 & 207 & 208 & 209 & 210 & 211 \\\hline 
167   & 5 & 2 & 2 & 1 & 1 & 1 & 2 & 0 & 2 & 2 & 2 & 1 & 0 & 1 & 2 & 2 & 2 & 2 & 2 \\
168   & 1 & 1 & 1 & 4 & 7 & 0 & 1 & 1 & 1 & 1 & 1 & 2 & 1 & 2 & 1 & 1 & 1 & 1 & 1 \\
169   & 1 & 1 & 1 & 4 & 5 & 5 & 1 & 1 & 1 & 1 & 1 & 2 & 1 & 2 & 1 & 1 & 1 & 1 & 1 \\
170   & 2 & 4 & 3 & 1 & 1 & 1 & 6 & 2 & 10 & 10 & 6 & 1 & 2 & 1 & 7 & 10 & 5 & 10 & 4 \\
171   & 1 & 1 & 1 & 9 & 4 & 4 & 1 & 1 & 1 & 1 & 1 & 2 & 1 & 2 & 1 & 1 & 1 & 1 & 1 \\
172   & 2 & 5 & 3 & 1 & 1 & 1 & 4 & 2 & 4 & 4 & 4 & 1 & 2 & 1 & 4 & 4 & 4 & 4 & 0 \\
173   & 1 & 1 & 1 & 4 & 6 & 6 & 1 & 1 & 1 & 1 & 1 & 2 & 1 & 2 & 1 & 1 & 1 & 1 & 1 \\
\end{tabular}
\vskip 20pt
{\large Table 20}
\end{center}
 }
\end{table}\clearpage

\begin{table}\tiny{
\begin{center}
\begin{tabular}{p{.35cm}|*{19}{p{.35cm}}}
    & 212 & 213 & 214 & 215 & 216 & 217 & 218 & 219 & 220 & 221 & 222 & 223 & 224 & 225 & 226 & 227 & 228 & 229 & 230 \\\hline 
1   & 1 & 1 & 1 & 1 & 1 & 1 & 1 & 1 & 2 & 1 & 1 & 0 & 6 & 3 & 2 & 9 & 6 & 6 & 5 \\
2   & 2 & 5 & 2 & 2 & 5 & 2 & 2 & 2 & 1 & 0 & 3 & 1 & 1 & 1 & 1 & 1 & 1 & 1 & 1 \\
3   & 2 & 7 & 2 & 2 & 0 & 2 & 2 & 2 & 1 & 5 & 3 & 1 & 1 & 1 & 1 & 1 & 1 & 1 & 1 \\
4   & 3 & 2 & 3 & 3 & 2 & 3 & 3 & 3 & 1 & 2 & 2 & 1 & 1 & 1 & 1 & 1 & 1 & 1 & 1 \\
5   & 1 & 1 & 1 & 1 & 1 & 1 & 1 & 1 & 2 & 1 & 1 & 4 & 4 & 3 & 2 & 4 & 4 & 4 & 4 \\
6   & 4 & 2 & 4 & 7 & 2 & 7 & 7 & 7 & 1 & 2 & 2 & 1 & 1 & 1 & 1 & 1 & 1 & 1 & 1 \\
7   & 5 & 2 & 5 & 4 & 2 & 4 & 4 & 4 & 1 & 2 & 2 & 1 & 1 & 1 & 1 & 1 & 1 & 1 & 1 \\
\end{tabular}
\vskip 20pt
{\large Table 21}
\end{center}
 }
\end{table}\clearpage

\begin{table}\tiny{
\begin{center}
\begin{tabular}{p{.35cm}|*{19}{p{.35cm}}}
    & 212 & 213 & 214 & 215 & 216 & 217 & 218 & 219 & 220 & 221 & 222 & 223 & 224 & 225 & 226 & 227 & 228 & 229 & 230 \\\hline 
84   & 2 & 2 & 2 & 2 & 2 & 2 & 2 & 2 & 1 & 2 & 2 & 1 & 1 & 1 & 1 & 1 & 1 & 1 & 1 \\
85   & 4 & 2 & 4 & 5 & 2 & 5 & 5 & 5 & 1 & 2 & 2 & 1 & 1 & 1 & 1 & 1 & 1 & 1 & 1 \\
86   & 4 & 2 & 4 & 6 & 2 & 6 & 6 & 6 & 1 & 2 & 2 & 1 & 1 & 1 & 1 & 1 & 1 & 1 & 1 \\
87   & 2 & 7 & 2 & 2 & 0 & 2 & 2 & 2 & 1 & 5 & 3 & 1 & 1 & 1 & 1 & 1 & 1 & 1 & 1 \\
88   & 1 & 1 & 1 & 1 & 1 & 1 & 1 & 1 & 2 & 1 & 1 & 6 & 6 & 3 & 2 & 6 & 6 & 9 & 5 \\
89   & 1 & 1 & 1 & 1 & 1 & 1 & 1 & 1 & 2 & 1 & 1 & 5 & 5 & 3 & 2 & 5 & 5 & 5 & 6 \\
90   & 4 & 2 & 4 & 6 & 2 & 6 & 6 & 6 & 1 & 2 & 2 & 1 & 1 & 1 & 1 & 1 & 1 & 1 & 1 \\
\end{tabular}
\vskip 20pt
{\large Table 22}
\end{center}
 }
\end{table}\clearpage

\begin{table}\tiny{
\begin{center}
\begin{tabular}{p{.35cm}|*{19}{p{.35cm}}}
    & 212 & 213 & 214 & 215 & 216 & 217 & 218 & 219 & 220 & 221 & 222 & 223 & 224 & 225 & 226 & 227 & 228 & 229 & 230 \\\hline 
167   & 2 & 7 & 2 & 2 & 0 & 2 & 2 & 2 & 1 & 5 & 3 & 1 & 1 & 1 & 1 & 1 & 1 & 1 & 1 \\
168   & 1 & 1 & 1 & 1 & 1 & 1 & 1 & 1 & 2 & 1 & 1 & 0 & 6 & 3 & 2 & 9 & 6 & 6 & 5 \\
169   & 1 & 1 & 1 & 1 & 1 & 1 & 1 & 1 & 2 & 1 & 1 & 5 & 5 & 3 & 2 & 5 & 5 & 5 & 6 \\
170   & 4 & 2 & 4 & 10 & 2 & 0 & 0 & 10 & 1 & 2 & 2 & 1 & 1 & 1 & 1 & 1 & 1 & 1 & 1 \\
171   & 1 & 1 & 1 & 1 & 1 & 1 & 1 & 1 & 2 & 1 & 1 & 4 & 4 & 3 & 2 & 4 & 4 & 4 & 4 \\
172   & 0 & 2 & 0 & 4 & 2 & 4 & 4 & 4 & 1 & 2 & 2 & 1 & 1 & 1 & 1 & 1 & 1 & 1 & 1 \\
173   & 1 & 1 & 1 & 1 & 1 & 1 & 1 & 1 & 2 & 1 & 1 & 6 & 6 & 3 & 2 & 6 & 6 & 9 & 5 \\
\end{tabular}
\vskip 20pt
{\large Table 23}
\end{center}
 }
\end{table}\clearpage

\begin{table}\tiny{
\begin{center}
\begin{tabular}{p{.35cm}|*{19}{p{.35cm}}}
    & 231 & 232 & 233 & 234 & 235 & 236 & 237 & 238 & 239 & 240 & 241 & 242 & 243 & 244 & 245 & 246 & 247 & 248 & 249 \\\hline 
1   & 4 & 1 & 1 & 1 & 1 & 1 & 1 & 1 & 1 & 1 & 1 & 2 & 4 & 3 & 1 & 1 & 6 & 2 & 6 \\
2   & 1 & 2 & 2 & 2 & 2 & 2 & 2 & 2 & 2 & 3 & 2 & 1 & 1 & 1 & 4 & 3 & 1 & 1 & 1 \\
3   & 1 & 2 & 2 & 2 & 2 & 2 & 2 & 2 & 2 & 3 & 2 & 1 & 1 & 1 & 4 & 3 & 1 & 1 & 1 \\
4   & 1 & 3 & 3 & 3 & 3 & 3 & 3 & 3 & 2 & 2 & 3 & 1 & 1 & 1 & 2 & 2 & 1 & 1 & 1 \\
5   & 6 & 1 & 1 & 1 & 1 & 1 & 1 & 1 & 1 & 1 & 1 & 2 & 5 & 3 & 1 & 1 & 4 & 2 & 4 \\
6   & 1 & 5 & 7 & 4 & 7 & 7 & 4 & 7 & 2 & 2 & 3 & 1 & 1 & 1 & 2 & 2 & 1 & 1 & 1 \\
7   & 1 & 4 & 4 & 5 & 4 & 4 & 5 & 4 & 2 & 2 & 3 & 1 & 1 & 1 & 2 & 2 & 1 & 1 & 1 \\
\end{tabular}
\vskip 20pt
{\large Table 24}
\end{center}
 }
\end{table}\clearpage

\begin{table}\tiny{
\begin{center}
\begin{tabular}{p{.35cm}|*{19}{p{.35cm}}}
    & 231 & 232 & 233 & 234 & 235 & 236 & 237 & 238 & 239 & 240 & 241 & 242 & 243 & 244 & 245 & 246 & 247 & 248 & 249 \\\hline 
84   & 1 & 2 & 2 & 2 & 2 & 2 & 2 & 2 & 4 & 2 & 2 & 1 & 1 & 1 & 2 & 2 & 1 & 1 & 1 \\
85   & 1 & 7 & 5 & 4 & 5 & 5 & 4 & 5 & 2 & 2 & 3 & 1 & 1 & 1 & 2 & 2 & 1 & 1 & 1 \\
86   & 1 & 5 & 6 & 4 & 6 & 6 & 4 & 6 & 2 & 2 & 3 & 1 & 1 & 1 & 2 & 2 & 1 & 1 & 1 \\
87   & 1 & 2 & 2 & 2 & 2 & 2 & 2 & 2 & 2 & 3 & 2 & 1 & 1 & 1 & 4 & 3 & 1 & 1 & 1 \\
88   & 4 & 1 & 1 & 1 & 1 & 1 & 1 & 1 & 1 & 1 & 1 & 2 & 4 & 3 & 1 & 1 & 6 & 2 & 9 \\
89   & 4 & 1 & 1 & 1 & 1 & 1 & 1 & 1 & 1 & 1 & 1 & 2 & 4 & 3 & 1 & 1 & 5 & 2 & 5 \\
90   & 1 & 5 & 6 & 4 & 6 & 6 & 4 & 6 & 2 & 2 & 3 & 1 & 1 & 1 & 2 & 2 & 1 & 1 & 1 \\
\end{tabular}
\vskip 20pt
{\large Table 25}
\end{center}
 }
 \end{table}\clearpage

\begin{table}\tiny{
\begin{center}
\begin{tabular}{p{.35cm}|*{19}{p{.35cm}}}
    & 231 & 232 & 233 & 234 & 235 & 236 & 237 & 238 & 239 & 240 & 241 & 242 & 243 & 244 & 245 & 246 & 247 & 248 & 249 \\\hline 
167   & 1 & 2 & 2 & 2 & 2 & 2 & 2 & 2 & 2 & 3 & 2 & 1 & 1 & 1 & 4 & 3 & 1 & 1 & 1 \\
168   & 4 & 1 & 1 & 1 & 1 & 1 & 1 & 1 & 1 & 1 & 1 & 2 & 4 & 3 & 1 & 1 & 6 & 2 & 6 \\
169   & 4 & 1 & 1 & 1 & 1 & 1 & 1 & 1 & 1 & 1 & 1 & 2 & 4 & 3 & 1 & 1 & 5 & 2 & 5 \\
170   & 1 & 5 & 8 & 4 & 0 & 10 & 4 & 10 & 2 & 2 & 3 & 1 & 1 & 1 & 2 & 2 & 1 & 1 & 1 \\
171   & 0 & 1 & 1 & 1 & 1 & 1 & 1 & 1 & 1 & 1 & 1 & 2 & 5 & 3 & 1 & 1 & 4 & 2 & 4 \\
172   & 1 & 4 & 4 & 0 & 4 & 4 & 8 & 4 & 2 & 2 & 3 & 1 & 1 & 1 & 2 & 2 & 1 & 1 & 1 \\
173   & 4 & 1 & 1 & 1 & 1 & 1 & 1 & 1 & 1 & 1 & 1 & 2 & 4 & 3 & 1 & 1 & 6 & 2 & 9 \\
\end{tabular}
\vskip 20pt
{\large Table 26}
\end{center}
 }
\end{table}

\clearpage

\section{Gauss codes} 

$ $

$3_{1}=\{1, -3, 5, -1, 3, -5\}$\newline

$4_{1}=\{-1, 3, -2, 4, -3, 1, -4, 2\}$\newline

$5_{1}=\{1, -3, 5, -7, 9, -1, 3, -5, 7, -9\}$\newline

$5_{2}=\{-1, 3, -5, 7, -9, 5, -3, 1, -7, 9\}$\newline

$6_{1}=\{-1, 3, -5, 7, -2, 4, -7, 5, -3, 1, -4, 2\}$\newline

$6_{2}=\{1, -3, 5, -7, 2, -4, 7, -1, 3, -5, 4, -2\}$\newline

$6_{3}=\{-2, 4, -1, 3, -5, 1, -6, 2, -4, 6, -3, 5\}$\newline

$7_{1}=\{ 1, -3, 5, -7, 9, -11, 13, -1, 3, -5, 7, -9, 11,
-13\}$\newline

$7_{2}=\{ -1, 3, -5, 7, -9, 13, -11, 9, -7, 5, -3, 1, -13,
11\}$\newline

$7_{3}=\{ -2, 4, -6, 8, -14, 12, -10, 2, -4, 6, -8, 10, -12,
14\}$\newline

$7_{4}=\{ 2, -4, 6, -8, 10, -12, 4, -2, 14, -10, 8, -6, 12,
-14\}$\newline

$7_{5}=\{ 1, -3, 5, -7, 9, -11, 13, -9, 7, -1, 3, -5, 11,
-13\}$\newline

$7_{6}=\{ 1, -3, 2, -4, 5, -7, 4, -2, 9, -1, 3, -9,  7, -5\}$\newline

$7_{7}=\{ 2, -4, 6, -8, 1, -3, 4, -2, 3, -5, 8, -6, 5, -1\}$\newline

$8_{1}=\{ -1, 3, -5, 7, -9, 11, -2, 4, -11, 9, -7, 5, -3, 1, -4,
2\}$\newline

$8_{2}=\{ 1, -3, 5, -7, 9, -11, 2, -4, 11, -1, 3, -5, 7, -9, 4,
-2\}$\newline

$8_{3}=\{ -1, 3, -5, 7, -2, 4, -6, 8, -7, 5, -3, 1, -8, 6, -4,
2\}$\newline

$8_{4}=\{ -2, 4, -6, 8, -1, 3, -5, 7, -8, 6, -4, 2, -3, 5, -7,
1\}$\newline

$8_{5}=\{ -2, 4, -6, 8, -10, 12, -1, 3, -8, 10, -12, 2, -4, 6, -3,
1\}$\newline

$8_{6}=\{ 1, -3, 5, -7, 9, -11, 2, -4, 11, -9, 7, -1, 3, -5, 4,
-2\}$\newline

$8_{7}=\{ -2, 4, -6, 8, -1, 3, -5, 1, -10, 2, -4, 6, -8, 10, -3,
5\}$\newline

$8_{8}=\{ 1, -3, 2, -4, 5, -1, 3, -5, 6, -8, 10, -2, 4, -10, 8,
-6\}$\newline

$8_{9}=\{ -2, 4, -6, 8, -1, 3, -5, 7, -8, 2, -4, 6, -7, 1, -3,
5\}$\newline

$8_{10}=\{ -2, 4, -6, 1, -3, 8, -10, 5, -1, 3, -5, 2, -4, 6, -8,
10\}$\newline

$8_{11}=\{ 1, -3, 2, -4, 5, -7, 9, -1, 3, -11, 4, -2, 11, -9, 7,
-5\}$\newline

$8_{12}=\{ -1, 3, -2, 4, -6, 8, -3, 1, -8, 6, -5, 7, -4, 2, -7,
5\}$\newline

$8_{13}=\{ 2, -4, 6, -8, 10, -2, 1, -3, 5, -1, 4, -10, 8, -6, 3,
-5\}$\newline

$8_{14}=\{ -1, 3, -5, 2, -4, 5, -7, 9, -11, 4, -2, 7, -3, 1, -9,
11\}$\newline

$8_{15}=\{ 1, -3, 5, -7, 9, -11, 13, -5, 7, -13, 15, -1, 3, -15, 11,
-9\}$\newline

$8_{16}=\{ 1, -3, 5, -7, 9, -2, 4, -5, 7, -6, 2, -1, 3, -4, 6,
-9\}$\newline

$8_{17}=\{ 1, -3, 5, -2, 4, -6, 8, -5, 7, -4, 6, -1, 3, -8, 2,
-7\}$\newline

$8_{18}=\{ -2, 1, -3, 4, -6, 5, -1, 8, -4, 7, -5, 2, -8, 3, -7,
6\}$\newline

$8_{19}=\{ -2, 4, -6, -8, 10, -12, 14, 6, -16, -10, 12, 2, -4, -14, 8,
16\}$\newline

$8_{20}=\{ 1, -3, -2, 5, 4, -6, -7, 2, 9, -4, 6, -1, 3, 7, -5,
-9\}$\newline

$8_{21}=\{ 1, -3, 5, 7, -9, 11, 2, -5, -4, 9, -11, -1, 3, -2, -7,
4\}$\newline

$9_{1}=\{ 1, -3, 5, -7, 9, -11, 13, -15, 17, -1, 3, -5, 7, -9, 11,
-13, 15, -17\}$\newline

$9_{2}=\{ -1, 3, -5, 7, -9, 11, -13, 15, -11, 9, -7, 5, -3, 1, -17,
13, -15, 17\}$\newline

$9_{3}=\{ -2, 4, -6, 8, -10, 12, -14, 16, -18, 2, -4, 6, -8, 10, -12,
18, -16, 14\}$\newline

$9_{4}=\{ -1, 3, -5, 7, -9, 11, -13, 15, -17, 9, -7, 5, -3, 1, -11,
13, -15, 17\}$\newline

$9_{5}=\{ -2, 4, -6, 8, -10, 12, -14, 16, -18, 10, -8, 6, -4, 2, -16,
14,-12, 18\}$\newline

$9_{6}=\{ 1, -3, 5, -7, 9, -11, 13, -15, 17, -13, 11, -1, 3, -5, 7,
-9, 15,-17\}$\newline

$9_{7}=\{1,-3,5,-7,9,-11,13,-15,17,-13,11,-9,7,-1,3,-5,15,-17\}$\newline

$9_{8 }= \{18,2,-4,6,-8,1,-3,8,-6,4,-2,5,-7,9,-5,3,-1,7,-9\}$\newline

$9_{9}= \{18, 1, -3, 5, -7, 9, -11, 13, -15, 17, -1, 3, -5, 7, -17,
15, -9, 11,-13\}$\newline

$9_{10}=\{ 2, -4, 6, -8, 10, -12, 14, -16, 18, -6, 4, -2, 8, -10, 12,
-18,16, -14\}$\newline

$9_{11}=\{ -2, 4, -6, 8, -10, 12, -14, 2, -1, 3, -8, 6, -3, 1, -4, 10,
-12,14\}$\newline

$9_{12}=\{ 1, -3, 2, -4, 5, -7, 9, -11, 4, -2, 13, -1, 3, -13, 11, -9,
7,-5\}$\newline

$9_{13}=\{ 2, -4, 6, -8, 10, -12, 14, -16, 18, -6, 4, -2, 8, -18, 16,
-10,12, -14\}$\newline

$9_{14}=\{ 2, -4, 6, -8, 1, -3, 10, -12, 3, -5, 8, -6, 4, -2, 5, -1,
12,-10\}$\newline

$9_{15}=\{-1, 2, -4, 1, -6, 8, -10, 12, -14, 10, -8, 6, -3, 4, -2, 3,
-12,14\}$\newline

$9_{16}=\{ 2, -4, 6, -8, 10, -12, 14, -16, 18, -2, 4, -14, 16, -18,
12, -6,8, -10\}$\newline

$9_{17}=\{ -1, 2, -4, 1, -3, 5, -7, 6, -8, 4, -2, 3, -5, 7, -9, 8, -6, 9\}$\newline

$9_{18}=\{ -1, 3, -5, 7, -9, 11, -13, 15, -17, 1, -15, 13, -7, 5, -3, 9,-11, 17\}$\newline

$9_{19}=\{ -1, 3, -5, 7, -2, 4, -6, 8, -7, 9, -4, 2, -9, 5, -3, 1, -8, 6\}$\newline

$9_{20}=\{ -1, 2, -4, 3, -5, 7, -3, 1, -9, 11, -13, 5, -7, 4, -2, 9, -11,13\}$\newline

$9_{21}=\{ -2, 4, -6, 8, -10, 12, -14, 6, -1, 3, -4, 2, -3, 1, -8, 14, -12,10\}$\newline

$9_{22}=\{ -1, 3, -5, 7, -2, 4, -6, 8, -3, 1, -4, 6, -8, 10, -7, 5, -10,2\}$\newline

$9_{23}=\{ -1, 3, -5, 7, -9, 11, -13, 1, -3, 9, -15, 17, -11, 13, -17, 15,-7, 5\}$\newline

$9_{24}=\{ 1, -3, 5, -1, 2, -7, 9, -2, 3, -5, 4, -6, 8, -9, 7, -4, 6, -8\}$\newline

$9_{25}=\{ 1, -3, 5, -7, 9, -1, 3, -9, 11, -13, 2, -4, 13, -11, 7, -5, 4,-2\}$\newline

$9_{26}=\{ -2, 4, -6, 8, -1, 3, -10, 12, -3, 5, -8, 2, -4, 6, -5, 1, -12,10\}$\newline

$9_{27}=\{ 1, -3, 2, -4, 6, -8, 5, -7, 8, -2, 9, -1, 3, -9, 4, -6, 7, -5\}$\newline

$9_{28}=\{ 1, -3, 2, -4, 5, -1, 3, -5, 6, -7, 9, -2, 4, -11, 7, -9, 11,-6\}$\newline

$9_{29}=\{ -1, 3, -5, 7, -9, 1, -2, 4, -7, 6, -8, 5, -4, 2, -3, 8, -6, 9\}$\newline

$9_{30}=\{ 1, -3, 2, -4, 6, -8, 3, -1, 5, -7, 4, -6, 9, -5, 7, -9, 8, -2\}$\newline

$9_{31}=\{ 1, -3, 5, -7, 2, -4, 9, -1, 3, -9, 6, -2, 11, -5, 7, -11, 4,-6\}$\newline

$9_{32}=\{ -2, 4, -6, 8, -1, 3, -10, 12, -3, 5, -8, 2, -12, 10, -4, 6, -5,1\}$\newline

$9_{33}=\{ 2, -1, 3, -2, 4, -6, 8, -5, 7, -4, 6, -9, 5, -3, 1, -7, 9, -8\}$\newline

$9_{34}=\{ -2, 1, -3, 4, -5, 7, -6, 3, -1, 8, -7, 9, -4, 2, -8, 6, -9, 5\}$\newline

$9_{35}=\{ -1, 3, -5, 7, -9, 11, -13, 15, -3, 1, -17, 9, -7, 5, -15, 13,-11, 17\}$\newline

$9_{36}=\{ -2, 4, -1, 3, -6, 8, -3, 1, -10, 12, -8, 6, -14, 2, -4, 10, -12,14\}$\newline

$9_{37}=\{ 2, -4, 6, -8, 1, -3, 5, -7, 4, -2, 7, -9, 8, -6, 9, -5, 3, -1\}$\newline

$9_{38}=\{ 1, -3, 5, -7, 9, -11, 3, -1, 13, -15, 17, -5, 11, -13, 15, -9,7, -17\}$\newline

$9_{39}=\{ 2, -4, 6, -8, 1, -3, 10, -12, 4, -2, 14, -10, 8, -6, 12, -14, 3,-1\}$\newline

$9_{40}=\{ 2, -1, 3, -4, 5, -7, 6, -3, 9, -2, 7, -11, 4, -9, 1, -6, 11,-5\}$\newline

$9_{41}=\{ -1, 3, -5, 7, -2, 4, -9, 11, -3, 1, -4, 6, -7, 5, -11, 9, -6,2\}$\newline

$9_{42}=\{ -2, -1, 3, 4, -6, 2, -5, 7, -8, 10, -7, 5, 1, -3, -10, 8, -4,6\}$\newline

$9_{43}=\{ 1, -3, 2, -4, 6, -8, 3, -1, -10, 12, 4, -6, -14, 10, -12, 14, 8,-2\}$\newline

$9_{44}=\{ 2, -4, 1, -3, -6, 5, 8, -1, 3, 7, -5, -9, 4, -2, 9, -8, -7, 6\}$\newline

$9_{45}=\{ 1, -3, 5, 7, -9, -1, 3, -5, -11, 13, -2, 9, -7, 4, -13, 11, -4,2\}$\newline

$9_{46}=\{ 2, -4, -1, 3, -5, 7, -9, 11, 4, -2, 6, 5, -3, 1, -11, 9, -7,-6\}$\newline

$9_{47}=\{ 1, -2, -4, 6, 8, -1, 3, -10, -6, 4, 12, -3, 5, -8, 2, -12, 10,-5\}$\newline

$9_{48}=\{ 2, -4, 6, -8, -10, 12, -1, -14, 4, -2, -3, 1, 8, -6, 14, 3, -12,10\}$\newline

$9_{49}=\{ 2, -4, 6, -8, -10, 12, 14, -16, 4, -2, 18, -14, 8, -6, 16, -18,-12, 10\}$\newline

$10_{1}=\{ -1, 3, -5, 7, -9, 11, -13, 2, -4, 13, -11, 9, -7, 5,
-3, 1, -15, 4, -2, 15\}$\newline

$10_{2}=\{ 1, -3, 5, -7, 9, -11, 13, -2, 4, -13, 15, -1, 3, -5,
7, -9, 11, -4, 2, -15\}$\newline

$10_{3}=\{ -1, 3, -5, 7, -9, 11, -2, 4, -6, 8, -11, 9, -7, 5,
-3, 1, -8, 6, -4, 2\}$\newline

$10_{4}=\{ 2, -4, 6, -8, 10, -12, 1, -3, 5, -7, 12, -10, 8, -6,
4, -2, 3, -5, 7, -1\}$\newline

$10_{5}=\{ -2, 4, -6, 8, -10, 12, -1, 3, -5, 1, -14, 2, -4, 6,
-8, 10, -12, 14, -3, 5\}$\newline

$10_{6}=\{ 1, -3, 5, -7, 9, -11, 13, -15, 2, -4, 15, -13, 11,
-1, 3, -5, 7, -9, 4, -2\}$\newline

$10_{7}=\{ -1, 3, -5, 7, -9, 11, -13, 15, -2, 4, -15, 9, -7, 5,
-3, 1, -11, 13, -4, 2\}$\newline

$10_{8}=\{ -1, 3, -5, 7, -9, 2, -4, 6, -8, 1, -3, 5, -7, 9,
-11, 8, -6, 4, -2, 11\}$\newline

$10_{9}=\{ -2, 4, -6, 8, -10, 12, -1, 3, -5, 7, -12, 2, -4, 6,
-8, 10, -7, 1, -3, 5\}$\newline

$10_{10}=\{ 2, -4, 6, -8, 10, -12, 14, -1, 3, -10, 8, -6, 4,
-2, 12, -5, 1, -3, 5, -14\}$\newline

$10_{11}=\{ 2, -4, 6, -8, 1, -3, 5, -7, 9, -11, 8, -6, 4, -2,
11, -9, 7, -1, 3, -5\}$\newline

$10_{12}=\{ -2, 4, -6, 8, -1, 3, -5, 1, -10, 12, -14, 2, -4, 6,
-8, 14, -12, 10, -3, 5\}$\newline

$10_{13}=\{ -1, 3, -5, 7, -9, 11, -2, 4, -11, 9, -6, 8, -7, 5,
-3, 1, -8, 6, -4, 2\}$\newline

$10_{14}=\{ 1, -3, 5, -7, 2, -4, 9, -11, 13, -1, 3, -5, 7, -13,
11, -15, 4, -2, 15, -9\}$\newline

$10_{15}=\{ -2, 4, -6, 8, -1, 3, -5, 7, -9, 5, -3, 1, -10, 2,
-4, 6, -8, 10, -7, 9\}$\newline

$10_{16}=\{ -1, 3, -5, 7, -2, 4, -6, 8, -10, 12, -7, 5, -3, 1,
-12, 2, -4, 10, -8, 6\}$\newline

$10_{17}=\{ -2, 4, -6, 8, -1, 3, -5, 7, -9, 1, -10, 2, -4, 6,
-8, 10, -3, 5, -7, 9\}$\newline

$10_{18}=\{ 1, -3, 2, -4, 6, -8, 3, -5, 7, -9, 11, -7, 5, -1,
8, -6, 4, -2, 9, -11\}$\newline

$10_{19}=\{ 2, -4, 1, -3, 5, -7, 9, -1, 6, -2, 8, -10, 4, -6,
3, -5, 7, -9, 10, -8\}$\newline

$10_{20}=\{ -1, 3, -5, 7, -9, 11, -13, 15, -2, 4, -11, 13, -15,
9, -7, 5, -3, 1, -4, 2\}$\newline

$10_{21}=\{ 1, -3, 5, -7, 9, -11, 13, -15, 2, -4, 15, -1, 3,
-5, 7, -13, 11, -9, 4, -2\}$\newline

$10_{22}=\{ -2, 4, -6, 8, -10, 12, -1, 3, -5, 7, -12, 10, -8,
2, -4, 6, -7, 1, -3, 5\}$\newline

$10_{23}=\{ -2, 4, -6, 8, -10, 12, -1, 3, -5, 1, -14, 6, -4, 2,
-8, 10, -12, 14, -3, 5\}$\newline

$10_{24}=\{ -2, 4, -1, 3, -5, 7, -9, 11, -13, 15, -4, 2, -11,
9, -7, 13, -15, 5, -3, 1\}$\newline

$10_{25}=\{ 1, -3, 5, -7, 9, -11, 2, -4, 11, -13, 15, -9, 7,
-1, 3, -5, 4, -2, 13, -15\}$\newline

$10_{26}=\{ 2, -4, 6, -8, 10, -1, 3, -5, 7, -6, 4, -2, 8, -10,
12, -3, 5, -7, 1, -12\}$\newline

$10_{27}=\{ -2, 4, -6, 2, -1, 3, -5, 7, -9, 11, -13, 1, -4, 6,
-3, 9, -11, 13, -7, 5\}$\newline

$10_{28}=\{ 2, -4, 6, -8, 10, -2, 12, -14, 1, -3, 5, -1, 4,
-10, 8, -6, 3, -5, 14, -12\}$\newline

$10_{29}=\{ 2, -1, 3, -5, 7, -2, 4, -6, 8, -3, 5, -7, 1, -9,
11, -8, 6, -11, 9, -4\}$\newline

$10_{30}=\{ -1, 3, -5, 7, -9, 11, -2, 4, -11, 13, -15, 5, -3,
1, -7, 15, -13, 9, -4, 2\}$\newline

$10_{31}=\{ -1, 3, -5, 2, -4, 6, -8, 12, -2, 7, -9, 8, -6, 4,
-12, 5, -3, 1, -7, 9\}$\newline

$10_{32}=\{ -1, 3, -5, 2, -4, 6, -8, 7, -9, 5, -11, 4, -6, 8,
-2, 11, -3, 1, -7, 9\}$\newline

$10_{33}=\{ -1, 3, -5, 7, -2, 4, -6, 8, -9, 5, -3, 1, -7, 9,
-10, 2, -8, 6, -4, 10\}$\newline

$10_{34}=\{ 2, -4, 6, -8, 10, -12, 14, -1, 3, -2, 4, -5, 1, -3,
5, -14, 12, -10, 8, -6\}$\newline

$10_{35}=\{ 1, -3, 2, -4, 6, -8, 5, -7, 10, -12, 7, -5, 3, -1,
12, -10, 8, -6, 4, -2\}$\newline

$10_{36}=\{ -1, 2, -4, 1, -3, 5, -7, 4, -2, 3, -9, 11, -13, 15,
-5, 7, -15, 13, -11, 9\}$\newline

$10_{37}=\{ 1, -3, 2, -4, 6, -8, 10, -6, 4, -2, 5, -7, 9, -1,
3, -9, 7, -5, 8, -10\}$\newline

$10_{38}=\{ -2, 4, -1, 3, -5, 7, -3, 1, -9, 11, -13, 15, -4, 2,
-15, 13, -11, 5, -7, 9\}$\newline

$10_{39}=\{ -1, 3, -2, 4, -5, 7, -9, 11, -13, 15, -4, 2, -15,
5, -7, 9, -3, 1, -11, 13\}$\newline

$10_{40}=\{ -2, 1, -3, 5, -1, 4, -6, 8, -10, 12, -14, 10, -8,
3, -5, 2, -4, 6, -12, 14\}$\newline

$10_{41}=\{ 1, -3, 2, -4, 5, -7, 6, -8, 7, -9, 11, -5, 3, -1,
8, -6, 9, -11, 4, -2\}$\newline

$10_{42}=\{ -1, 3, -2, 5, -7, 4, -6, 8, -10, 2, -4, 9, -5, 7,
-9, 10, -3, 1, -8, 6\}$\newline

$10_{43}=\{ 1, -3, 2, -4, 6, -8, 5, -7, 4, -2, 9, -1, 3, -9, 7,
-5, 10, -6, 8, -10\}$\newline

$10_{44}=\{ 1, -3, 2, -4, 5, -7, 9, -1, 3, -9, 6, -8, 7, -11,
4, -2, 11, -5, 8, -6\}$\newline

$10_{45}=\{ 2, -1, 3, -2, 4, -5, 7, -6, 8, -3, 1, -4, 10, -7,
9, -8, 6, -9, 5, -10\}$\newline

$10_{46}=\{ -2, 4, -6, 8, -10, 12, -14, 16, -1, 3, -12, 14,
-16, 2, -4, 6, -8, 10, -3, 1\}$\newline

$10_{47}=\{ -2, 4, -6, 8, -10, 1, -3, 12, -14, 5, -1, 3, -5, 2,
-4, 6, -8, 10, -12, 14\}$\newline

$10_{48}=\{ -1, 3, -5, 7, -2, 4, -6, 8, -10, 2, -4, 6, -9, 1,
-3, 5, -7, 9, -8, 10\}$\newline

$10_{49}=\{ -1, 3, -5, 7, -9, 11, -13, 15, -17, 9, -11, 17,
-19, 1, -3, 5, -7, 19, -15, 13\}$\newline

$10_{50}=\{ -2, 4, -6, 8, -10, 12, -14, 16, -1, 3, -12, 14,
-16, 6, -4, 2, -8, 10, -3, 1\}$\newline

$10_{51}=\{ -2, 4, -6, 8, -10, 1, -3, 5, -1, 12, -14, 3, -5, 6,
-4, 2, -8, 10, -12, 14\}$\newline

$10_{52}=\{ -1, 3, -5, 7, -2, 4, -6, 8, -10, 2, -4, 6, -9, 5,
-3, 1, -7, 9, -8, 10\}$\newline

$10_{53}=\{ -1, 3, -5, 7, -9, 11, -13, 15, -17, 9, -11, 17,
-19, 5, -3, 1, -7, 19, -15, 13\}$\newline

$10_{54}=\{ -1, 3, -2, 4, -6, 8, -10, 2, -4, 6, -5, 7, -9, 1,
-3, 9, -7, 5, -8, 10\}$\newline

$10_{55}=\{ -1, 3, -5, 7, -9, 11, -13, 5, -7, 13, -15, 17, -19,
1, -3, 19, -17, 15, -11, 9\}$\newline

$10_{56}=\{ -2, 4, -6, 8, -10, 12, -1, 3, -8, 10, -12, 14, -16,
2, -4, 16, -14, 6, -3, 1\}$\newline

$10_{57}=\{ -2, 4, -6, 1, -3, 5, -1, 8, -10, 3, -5, 12, -14, 2,
-4, 14, -12, 6, -8, 10\}$\newline

$10_{58}=\{ -1, 3, -5, 7, -2, 4, -6, 8, -7, 5, -9, 11, -8, 6,
-11, 9, -3, 1, -4, 2\}$\newline

$10_{59}=\{ -1, 3, -2, 4, -3, 1, -6, 8, -5, 7, -8, 10, -12, 6,
-7, 5, -4, 2, -10, 12\}$\newline

$10_{60}=\{ -2, 1, -3, 5, -7, 2, -4, 3, -1, 4, -6, 9, -11, 7,
-5, 6, -8, 11, -9, 8\}$\newline

$10_{61}=\{ -1, 3, -5, 2, -4, 6, -8, 10, -12, 5, -3, 1, -7, 8,
-10, 12, -2, 4, -6, 7\}$\newline

$10_{62}=\{ -2, 4, -6, 8, -10, 12, -1, 3, -14, 2, -4, 6, -5, 1,
-3, 5, -8, 10, -12, 14\}$\newline

$10_{63}=\{ -1, 3, -5, 7, -9, 11, -7, 13, -15, 17, -13, 5, -3,
1, -19, 15, -17, 9, -11, 19\}$\newline

$10_{64}=\{ -2, 4, -1, 3, -5, 7, -6, 2, -4, 8, -10, 12, -7, 1,
-3, 5, -8, 10, -12, 6\}$\newline

$10_{65}=\{ -2, 4, -6, 8, -1, 3, -10, 12, -14, 6, -4, 2, -8,
10, -12, 14, -5, 1, -3, 5\}$\newline

$10_{66}=\{ -1, 3, -5, 7, -9, 11, -13, 9, -15, 17, -19, 15, -7,
1, -3, 5, -17, 19, -11, 13\}$\newline

$10_{67}=\{ -1, 3, -5, 7, -9, 11, -2, 4, -11, 13, -15, 1, -3,
15, -13, 9, -7, 5, -4, 2\}$\newline

$10_{68}=\{ -2, 4, -1, 3, -5, 7, -9, 11, -13, 1, -6, 2, -4, 6,
-7, 5, -3, 13, -11, 9\}$\newline

$10_{69}=\{ -2, 4, -1, 3, -6, 8, -3, 5, -10, 12, -5, 1, -14, 2,
-4, 14, -8, 6, -12, 10\}$\newline

$10_{70}=\{ -2, 4, -6, 8, -10, 1, -3, 6, -8, 10, -5, 7, -4, 2,
-7, 5, -12, 3, -1, 12\}$\newline

$10_{71}=\{ -2, 4, -1, 3, -6, 8, -5, 1, -3, 5, -7, 9, -4, 2,
-9, 7, -10, 6, -8, 10\}$\newline

$10_{72}=\{ 2, -4, 6, -8, 10, -12, 14, -6, 8, -10, 1, -3, 4,
-16, 12, -14, 16, -2, 3, -1\}$\newline

$10_{73}=\{ -1, 3, -2, 5, -7, 9, -5, 11, -13, 7, -9, 4, -3, 1,
-4, 2, -6, 13, -11, 6\}$\newline

$10_{74}=\{ -1, 3, -5, 2, -4, 7, -9, 11, -3, 1, -13, 5, -15, 4,
-2, 15, -11, 9, -7, 13\}$\newline

$10_{75}=\{ 1, -2, 4, -1, 3, -6, 8, -10, 12, -4, 2, -3, 5, -8,
6, -5, 7, -12, 10, -7\}$\newline

$10_{76}=\{ -2, 4, -6, 8, -1, 3, -8, 6, -10, 2, -4, 12, -14,
16, -3, 1, -12, 14, -16, 10\}$\newline

$10_{77}=\{ -2, 4, -6, 8, -4, 2, -10, 12, -14, 1, -3, 5, -1, 6,
-8, 3, -5, 10, -12, 14\}$\newline

$10_{78}=\{ 1, -3, 5, -7, 9, -1, 3, -9, 2, -4, 11, -13, 15,
-11, 7, -5, 13, -15, 4, -2\}$\newline

$10_{79}=\{ -2, 4, -6, 1, -3, 5, -7, 9, -1, 3, -5, 8, -10, 7,
-9, 2, -4, 6, -8, 10\}$\newline

$10_{80}=\{ -1, 3, -5, 7, -9, 11, -13, 5, -7, 9, -15, 17, -3,
1, -19, 15, -17, 19, -11, 13\}$\newline

$10_{81}=\{ -2, 4, -6, 8, -10, 2, -4, 10, -1, 3, -5, 7, -9, 1,
-3, 9, -8, 6, -7, 5\}$\newline

$10_{82}=\{ 1, -3, 5, -7, 9, -2, 4, -6, 8, -9, 11, -4, 6, -1,
3, -5, 7, -8, 2, -11\}$\newline

$10_{83}=\{ 2, -4, 6, -8, 10, -12, 14, -2, 1, -3, 8, -14, 12,
-10, 5, -1, 4, -6, 3, -5\}$\newline

$10_{84}=\{ 2, -4, 1, -3, 6, -8, 5, -1, 10, -6, 8, -12, 14,
-10, 3, -5, 4, -2, 12, -14\}$\newline

$10_{85}=\{ 2, -1, 3, -4, 6, -5, 7, -9, 11, -2, 4, -13, 5, -7,
9, -11, 1, -3, 13, -6\}$\newline

$10_{86}=\{ 1, -3, 5, -7, 2, -4, 6, -8, 10, -12, 7, -1, 8, -6,
4, -10, 3, -5, 12, -2\}$\newline

$10_{87}=\{ 2, -4, 6, -8, 4, -2, 1, -3, 10, -12, 5, -1, 3, -7,
12, -6, 8, -10, 7, -5\}$\newline

$10_{88}=\{ 1, -3, 2, -4, 3, -5, 6, -7, 9, -8, 5, -1, 10, -9,
7, -10, 4, -2, 8, -6\}$\newline

$10_{89}=\{ -1, 3, -5, 7, -2, 9, -11, 4, -7, 5, -4, 13, -9, 6,
-3, 1, -6, 11, -13, 2\}$\newline

$10_{90}=\{ -2, 4, -6, 1, -3, 2, -4, 5, -7, 6, -8, 10, -12, 3,
-5, 7, -1, 12, -10, 8\}$\newline

$10_{91}=\{ 1, -3, 5, -2, 4, -6, 8, -10, 2, -7, 9, -8, 10, -1,
3, -5, 7, -4, 6, -9\}$\newline

$10_{92}=\{ 2, -4, 6, -8, 10, -12, 8, -2, 1, -3, 14, -10, 12,
-16, 3, -1, 4, -6, 16, -14\}$\newline

$10_{93}=\{ 1, -3, 5, -7, 9, -1, 2, -4, 6, -8, 7, -9, 10, -6,
4, -2, 3, -5, 8, -10\}$\newline

$10_{94}=\{ 2, -4, 6, -8, 10, -12, 1, -3, 5, -7, 8, -10, 12,
-2, 3, -5, 4, -6, 7, -1\}$\newline

$10_{95}=\{ 1, -3, 2, -4, 6, -8, 5, -1, 3, -5, 10, -12, 14, -2,
8, -10, 12, -6, 4, -14\}$\newline

$10_{96}=\{ 2, -4, 1, -6, 8, -3, 4, -2, 10, -12, 3, -1, 5, -7,
12, -10, 7, -8, 6, -5\}$\newline

$10_{97}=\{ 2, -4, 1, -3, 4, -6, 8, -10, 12, -14, 16, -8, 6,
-2, 14, -12, 3, -1, 10, -16\}$\newline

$10_{98}=\{ 1, -3, 2, -4, 5, -7, 9, -1, 3, -11, 13, -15, 4, -2,
11, -9, 7, -13, 15, -5\}$\newline

$10_{99}=\{ 1, -3, 5, -7, 9, -2, 4, -1, 3, -6, 8, -10, 2, -4,
6, -5, 7, -8, 10, -9\}$\newline

$10_{100}=\{ 1, -2, 4, -3, 5, -6, 2, -7, 9, -11, 3, -5, 13, -1,
7, -9, 11, -4, 6, -13\}$\newline

$10_{101}=\{ -2, 4, -6, 8, -10, 12, -14, 6, -16, 18, -8, 14,
-20, 2, -4, 20, -12, 10, -18, 16\}$\newline

$10_{102}=\{ -2, 4, -6, 1, -3, 8, -10, 5, -7, 6, -4, 2, -12, 3,
-5, 7, -1, 12, -8, 10\}$\newline

$10_{103}=\{ -1, 3, -2, 4, -5, 7, -6, 2, -9, 11, -4, 6, -3, 1,
-13, 9, -11, 5, -7, 13\}$\newline

$10_{104}=\{ -2, 1, -3, 5, -7, 4, -6, 9, -1, 8, -10, 3, -5, 7,
-9, 2, -8, 10, -4, 6\}$\newline

$10_{105}=\{ -2, 4, -6, 8, -1, 3, -5, 7, -8, 10, -12, 2, -4,
12, -3, 1, -10, 6, -7, 5\}$\newline

$10_{106}=\{ -2, 4, -6, 1, -3, 5, -7, 6, -8, 10, -12, 2, -4, 7,
-1, 8, -10, 3, -5, 12\}$\newline

$10_{107}=\{ -1, 3, -2, 4, -5, 7, -6, 8, -9, 1, -3, 9, -10, 6,
-4, 2, -8, 10, -7, 5\}$\newline

$10_{108}=\{ -2, 4, -1, 3, -5, 7, -6, 2, -4, 8, -10, 6, -9, 5,
-3, 1, -8, 10, -7, 9\}$\newline

$10_{109}=\{ -1, 3, -5, 7, -9, 2, -4, 6, -8, 10, -2, 1, -3, 4,
-6, 5, -7, 8, -10, 9\}$\newline

$10_{110}=\{ -2, 4, -6, 8, -1, 3, -5, 7, -9, 11, -8, 6, -7, 5,
-4, 2, -3, 9, -11, 1\}$\newline

$10_{111}=\{ -2, 4, -6, 8, -10, 12, -4, 2, -14, 1, -3, 10, -12,
14, -16, 6, -8, 3, -1, 16\}$\newline

$10_{112}=\{ -1, 3, -5, 2, -4, 7, -9, 6, -2, 11, -7, 8, -6, 1,
-3, 5, -11, 4, -8, 9\}$\newline

$10_{113}=\{ 2, -4, 6, -8, 1, -3, 4, -10, 12, -2, 3, -5, 8,
-14, 10, -12, 14, -6, 5, -1\}$\newline

$10_{114}=\{ -2, 1, -3, 4, -6, 5, -1, 8, -4, 7, -9, 11, -5, 2,
-8, 3, -11, 9, -7, 6\}$\newline

$10_{115}=\{ -1, 2, -4, 6, -8, 10, -2, 3, -5, 7, -9, 1, -3, 4,
-10, 9, -7, 8, -6, 5\}$\newline

$10_{116}=\{ 1, -3, 5, -7, 9, -2, 4, -1, 3, -6, 8, -9, 11, -4,
6, -5, 7, -8, 2, -11\}$\newline

$10_{117}=\{ -2, 4, -6, 8, -10, 12, -14, 1, -3, 10, -4, 2, -12,
5, -1, 6, -8, 3, -5, 14\}$\newline

$10_{118}=\{ -1, 3, -5, 2, -4, 7, -9, 6, -8, 5, -7, 10, -6, 1,
-3, 8, -2, 4, -10, 9\}$\newline

$10_{119}=\{ -2, 4, -6, 1, -3, 8, -4, 2, -10, 5, -7, 6, -8, 10,
-12, 7, -1, 3, -5, 12\}$\newline

$10_{120}=\{ -1, 3, -5, 7, -9, 11, -13, 5, -3, 15, -11, 17,
-19, 13, -15, 1, -7, 19, -17, 9\}$\newline

$10_{121}=\{ -1, 3, -5, 7, -2, 4, -3, 9, -11, 5, -6, 2, -13,
11, -9, 1, -4, 6, -7, 13\}$\newline

$10_{122}=\{ 1, -2, 4, -3, 5, -6, 8, -10, 2, -7, 3, -14, 6, -1,
7, -4, 10, -8, 14, -5\}$\newline

$10_{123}=\{ 1, -2, 4, -5, 7, -6, 8, -1, 9, -4, 10, -7, 11, -8,
2, -9, 5, -10, 6, -11\}$\newline

$10_{124}=\{ -2, 4, -6, 8, -10, 12, -14, 16, 18, -20, -12, 14,
-16, 2, -4, 6, -8, 10, 20, -18\}$\newline

$10_{125}=\{ -2, 4, -6, 8, -10, -1, 3, -5, -7, 9, 1, -3, 5, 2,
-4, 6, -8, 10, -9, 7\}$\newline

$10_{126}=\{ 1, -3, 5, -7, 9, 2, -4, 6, -11, 13, -2, 4, -6, -1,
3, -5, 7, -9, -13, 11\}$\newline

$10_{127}=\{ -1, 3, -5, 7, -3, -2, 4, 5, -7, 9, -11,13, -15,
-4,2, 1, -9, 11, -13, 15\}$\newline

$10_{128}=\{ -2, 4, -6, 8, -10, 12, -14, -16, 18, 6, -8, -20,
16,2, -4, -18, 20, 14, -12, 10\}$\newline

$10_{129}=\{ 1, -3, -2, 4, -6, 8, -10, 5, -7, 2, -4, 9, -5, -1,
3,7, -9, 10, -8, 6\}$\newline

$10_{130}=\{ -1, 3, -5, 7, -2, 4, -6, -9, 11, 2, -4, 6, -13, 5,
-3, 1, -7, 13, 9, -11\}$\newline

$10_{131}=\{ -1, 3, -5, 7, 9, -11, 13, 2, -4, -9, 11, -13, -15,
5,-3, 1, -7, 15, -2, 4\}$\newline

$10_{132}=\{ -1, -2, 4, 3, -6, -5, 2, 7, -9, 1, 5, -11, -3, 13,
-7, 9, -13, -4, 11, 6\}$\newline

$10_{133}=\{ 2, 1, -3, -4, 5, -7, -1, 9, -11, -2, 7, -13, 4,
15,-9, 11, -15, 3, 13, -5\}$\newline

$10_{134}=\{ 2, -4, 6, -8, -10, 12, 4, -14, 16, -2, -12, 18, 8,
-20, 14, -16, 20, -6, -18, 10\}$\newline

$10_{135}=\{ -2, 4, -6, 1, -3, 5, -1, -7, 9, 3, -5, 8, -10, 2,
-4, 10, -8, 6,7, -9\}$\newline

$10_{136}=\{ -1, 3, -5, 7, -2, 4, -3, -6, 8, 1, -4, 10, -7,
-12,6, -8, 12, 5, -10, 2\}$\newline

$10_{137}=\{ -2, 1, -3, 5, -7, 2, -4, 3, -1, 4, 9, -6, 8, 7,
-5,-9, 11, -8, 6, -11\}$\newline

$10_{138}=\{ -2, 1, -3, -4, 6, 2, -8, 3, -1, 8, -10, 5, -7, -6,
4,10, -12, 7, -5, 12\}$\newline

$10_{139}=\{ -2, 4, -6, -8, 10, -12, 14, 2, -4, -16, 8, 18,
-20,-10, 12, -14, 16, 6, -18, 20\}$\newline

$10_{140}=\{ 1, 2, -4, 6, -3, 5, -7, 9, -2, 4, -6, -11, 13, -1,
-9, 7, -5, 3, 11, -13\}$\newline

$10_{141}=\{ -2, 1, 4, -3, 5, -6, -7, 2, 9, -4, 8, -5, 11, 7,
-1,-9, 3, -8, 6, -11\}$\newline

$10_{142}=\{ 2, 4, -6, 8, -10, 12, -14, -2, 16, -18, 20, 10,
-12,14, -4, 6, -8, -20, 18, -16\}$\newline

$10_{143}=\{ 1, -3, -2, 4, -6, 5, -7, -1, 3, 9, -11, 13, -5, 7,
-9, 2, -4, 11, -13, 6\}$\newline

$10_{144}=\{ -1, 3, 2, -4, 6, -8, -5, 1, -3, 5, -7, 9, 8, -6,
4,-2, -11, 7, -9, 11\}$\newline

$10_{145}=\{ 1, 3, -5, -7, 9, 11, -13, -15, 17, -9, 7, -1, 15,
-2,-11, 5, -3, 13, 2, -17\}$\newline

$10_{146}=\{ -1, 3, 2, -4, 6, 8, -10, 5, -7, -6, 4, -2, -9, 10,
-8, 9, -3, 1, -5, 7\}$\newline

$10_{147}=\{ -1, 3, 2, -4, 5, 6, -3, -8, 10, 1, -6, -7, 4, -12,
8, -10, 12, -2, 7, -5\}$\newline

$10_{148}=\{ -2, 4, 1, -3, -6, 5, -7, 9, -11, -1, 3, 13, -5, 2,
-4, 7, -9, 11, -13, 6\}$\newline

$10_{149}=\{ 1, -3, 5, 7, -9, 11, -13, 15, -7, 9, -11, 2, -4,
13, -15, -1, 3, -5, -2, 4\}$\newline

$10_{150}=\{ -1, 3, 2, -4, -6, 8, -3, -10, 12, 1, -8, 14, 4,
-16, 10, -12, 16, -2, -14, 6\}$\newline

$10_{151}=\{ -2, 4, -6, 8, -10, 2, -4, 10, -1, 3, 12, -14, -5,
1, -3, 5, -8, 6, 14, -12\}$\newline

$10_{152}=\{ 1, -3, 5, 7, -9, 11, -13, 15, -7, 9, -11, -17, 19,
13, -15, -1, 3, -5, 17, -19\}$\newline

$10_{153}=\{ -1, 3, -5, 7, -9, 1, -3, 5, 2, -4, 6, -8, 10, -2,
4, -10, -7, 9, 8, -6\}$\newline

$10_{154}=\{ -2, 4, -6, 8, -10, 2, -4, 10, 12, -14, 16, -18,
20, -12, 14, -20, -8, 6, 18, -16\}$\newline

$10_{155}=\{ 1, -3, 5, -7, -2, 4, -6, -8, 10, -5, 7, -12, 8,
-1, 3, -10, 12, 2, -4, 6\}$\newline

$10_{156}=\{ 2, -4, 6, 1, -3, -5, 7, -9, 11, 3, -13, -2, 4, -7,
5, 13, -1, -11, 9, -6\}$\newline

$10_{157}=\{ -2, 4, -6, 8, 1, -10, 12, 2, -4, 3, 14, -16, 10,
-12, -3, 6, -8, -14, 16, -1\}$\newline

$10_{158}=\{ 1, -3, 5, -7, -2, 4, -6, -8, 10, -5, 7, -12, 8,
-1, 3, -10, 12, 6, -4, 2\}$\newline

$10_{159}=\{ 1, 3, -5, 7, 2, -9, 11, -4, -3, 5, -7, -13, 9, -6,
4, -1, 13, -2, 6, -11\}$\newline

$10_{160}=\{ -2, 1, -3, 4, -6, -8, 10, -12, -1, 14, -4, -16, 8,
2, -14, 3, 12, -10, 16, 6\}$\newline

$10_{161}=\{ 1, -3, 5, -7, -2, 9, -11, 13, -15, -5, 7, 17, -9,
-1, 3, 11, -13, 15, -17, 2\}$\newline

$10_{163}=\{ 2, -4, 6, 1, -8, -3, 5, -7, 9, 8, -11, -6, 4, -2,
7, -9, -1, 11, 3, -5\}$\newline

$10_{164}=\{ -2, 1, -3, 4, -6, 5, -1, 8, -4, -10, 12, -14, -5,
2, -8, 3, 14, -12, 10, 6\}$\newline

$10_{165}=\{ 1, 3, -5, -7, 9, -2, 4, -6, -3, 8, 7, -10, 2, -1,
-8, 5, 6, -4, 10, -9\}$\newline

$10_{166}=\{ 2, 1, -3, -4, 6, -8, 10, 3, -1, -12, 8, -14, 16,
-10, 12, -2, 4, -16, 14, -6\}$\newline


\begin{thebibliography}{99}


\bibitem{Birman}
        J. Birman, \emph{Braids, Links, and Mapping Class Groups}, Annals of Math. Studies, no. 82, Princeton University Press, 1974


\bibitem{CetS}
        J. S. Carter, M. Saito, \emph{Knotted Surfaces and Their Diagrams}, Mathematical Surveys and Monographs, vol. 55, American Mathematical Society, 1998

\bibitem{jsCetal}
        J. S. Carter, D. Jelsovsky, S. Kamada, L. Langford, M. Saito, \emph{Quandle Cohomology and State-sum Invariants of Knotted Curves and Surfaces}, math.GT/9903135, 1999
 
\bibitem{jsCetal1}
        J. Scott Carter, D. Jelsovsky, S. Kamada, M. Saito, \emph{Computations of Quandle Cocycle Invariants of Knotted Curves and Surfaces}, Adv. in Math. 157, 36-94, 2001


\bibitem{FR}
        R. Fenn, C. Rourke, \emph{Racks and Links in Codimension Two}, J. Knot Theory Ramifications, {\bf1}, (1992), no. 4, 343-406


\bibitem{J}
        V. F. R. Jones, \emph{Hecke Algebra Representations of Braid Groups and Link Polynomials}, Ann. of Math., {\bf126}, (1982), no. 2, 335-388


\bibitem{dJoyce}
        D. Joyce, \emph{A Classifying Invariant of Knots, The Knot Quandle}, J. Pure Appl. Alg., {\bf23}, (1982), 37-65

\bibitem{lhKauffman}
        L. H. Kauffman, \emph{Knots and Physics}, 2nd edition, Series on Knots and Everything, vol. 1, World Scientific, 1993

\bibitem{pLopes}
        P. Lopes, \emph{Quandles at Finite Temperatures I}, math.QA/0105099, 2001, accepted for publication in J. Knot Theory Ramifications.

\bibitem{sMatveev}
        S. V. Matveev, \emph{Distributive Groupoids in Knot Theory}, Math. USSR Sbornik, {\bf47}, (1984), no. 1, 73-83


\bibitem{R}
        D. Rolfsen, \emph{Knots and Links}, Mathematics Lecture Series, vol. 7, Publish or Perish, Inc., Houston, Texas, 1976


\bibitem{dRoseman1}
        D. Roseman, \emph{Reidemeister-Type Moves for Surfaces in Four-Dimensional Space}, in ``Knot Theory'', (V. F. R. Jones et al, Eds.), Banach Center Publications, vol.42, pp.347-380, Polish Academy of Sciences, Warsaw, 1998

\end{thebibliography}
\end{document}